\documentclass[12pt]{article}
\usepackage[square,comma,numbers]{natbib}
\usepackage{amsmath,amsfonts,amsthm,amssymb}
\usepackage{mathtools}
\usepackage{url,doi}
\usepackage[shortlabels]{enumitem}

\newcommand{\F}{\mathbb{F}}
\newcommand{\N}{\mathbb{N}}
\newcommand{\R}{\mathbb{R}}
\newcommand{\Z}{\mathbb{Z}}
\newcommand{\tabs}[1]{\mathopen|#1\mathclose|}
\newcommand{\restrict}{\upharpoonright}
\newcommand{\ssm}{\smallsetminus}
\newcommand{\partialto}{\rightharpoonup}
\newcommand{\concat}{{}^{\smallfrown}}

\DeclareMathOperator{\dom}{dom}

\DeclareMathOperator{\diam}{diam}
\DeclareMathOperator{\Sch}{Sch}

\newtheorem{theorem}{Theorem}[section]
\newtheorem{lemma}[theorem]{Lemma}
\newtheorem{corollary}[theorem]{Corollary}
\newtheorem{question}[theorem]{Question}

\theoremstyle{definition}
\newtheorem{definition}[theorem]{Definition}

\numberwithin{equation}{section}

\allowdisplaybreaks

\begin{document}
\title{Measurable domatic partitions}
\author{Edward Hou\\
\small Department of Mathematics\\
\small California Institute of Technology\\
\small \href{mailto:ehou@caltech.edu}{\nolinkurl{ehou@caltech.edu}}}
\date{}
\maketitle

\begin{abstract}
Let $\Gamma$ be a compact Polish group of finite topological dimension. For a countably infinite subset $S\subseteq \Gamma$, a \emph{domatic $\aleph_0$-partition} (for its Schreier graph on $\Gamma$) is a partial function $f:\Gamma\rightharpoonup\mathbb{N}$ such that for every $x\in \Gamma$, one has $f[S\cdot x]=\mathbb{N}$. We show that a continuous domatic $\aleph_0$-partition exists, if and only if a Baire measurable domatic $\aleph_0$-partition exists, if and only if the topological closure of $S$ is uncountable. A Haar measurable domatic $\aleph_0$-partition exists for all choices of $S$. We also investigate domatic partitions in the general descriptive graph combinatorial setting. 
\end{abstract}

\section{Introduction}

This work is concerned with the existence of domatic partitions in the area of descriptive graph combinatorics. 

Let $G$ be a directed graph on a vertex set $V$ with possible loop edges, and we represent its edge set as a binary relation $G\subseteq V^2$. For a vertex $v\in V$, its out-neighborhood is $N_G(v)=\{w\in V:(v,w)\in G\}$. A \emph{domatic partition} for $G$ is a partial function $f:V\partialto C$ which colors the vertices such that for every vertex $v\in V$, its out-neighborhood $N_G(v)$ is fully colored by $f$, meaning $f[N_G(v)]=C$. 

A classic result of Zelinka \citep{zelinka} can be stated as follows: Let $Q_n$ be the finite hypercube graph on $2^n$ vertices. Assume $Q_n$ is loop-free simple undirected, so that it is $n$-regular. Then $Q_n$ admits a domatic $n$-partition $f:V(Q_n)\to\{0,1,\ldots,n-1\}$ if and only if $n$ is a power of two. We are motivated by this result to find analogous criteria on the existence of domatic partitions for infinite graphs in the context of descriptive graph combinatorics. 

The area of descriptive graph combinatorics studies measurable combinatorial objects: Typically one defines a graph with a Borel measurable edge binary relation on a Polish space of vertices, and asks for the existence of special coloring functions that are measurable in certain senses (Borel, Baire, or measure). See for example the survey by Kechris--Marks \citep{kechris-marks}. 

In Section \ref{sec:2} we will be analyzing the Schreier graphs $\Sch(\Gamma,S,\Gamma)$: Given a Polish group $\Gamma$ and an arbitrary subset $S\subseteq\Gamma$, the \emph{Schreier graph} $G=\Sch(\Gamma,S,\Gamma)$ is the directed graph on the vertex set $V=\Gamma$ defined by the edge set $G=\{(\gamma,s\cdot\gamma):\gamma\in\Gamma,s\in S\}$. A domatic $\aleph_0$-partition in this case is a partial function $f:\Gamma\partialto{}\N$ such that $f[S\cdot\gamma]=\N$ for every $\gamma\in\Gamma$. We will prove the following results concerning the existence of various measurable kinds of domatic $\aleph_0$-partitions for $\Sch(\Gamma,S,\Gamma)$:

\begin{theorem}[Corollary \ref{cor:2.18}]\label{thm:1.1}
Let $\Gamma$ be a finite-dimensional compact Polish group, and let $S\subseteq\Gamma$ be a subset. Then the graph $\Sch(\Gamma,S,\Gamma)$ admits a domatic $\aleph_0$-partition with open parts, if and only if it admits a domatic $\aleph_0$-partition with Baire measurable parts, if and only if $\overline{S}\subseteq\Gamma$ is uncountable. 
\end{theorem}

The phrase ``finite-dimensional'' here means that the Polish group $\Gamma$ has finite Lebesgue covering dimension as a Polish space. In the case of compact Polish groups, we will give an alternative characterization of the dimension of $\Gamma$ in Definition \ref{def:2.4}. We leave open Question \ref{que:2.20} on the existence of domatic $\aleph_0$-partitions in infinite-dimensional compact Polish groups. 

\begin{theorem}[Corollary \ref{cor:2.19}]\label{thm:1.2}
Let $n\in\N$, and let $S\subseteq\R^n$ be a subset of $\R^n$. Then the graph $\Sch(\R^n,S,\R^n)$ admits a domatic $\aleph_0$-partition with open or Baire measurable parts if and only if either $\overline{S}\subseteq\R^n$ is uncountable or $S\subseteq\R^n$ is unbounded. 
\end{theorem}

\begin{theorem}[Corollary \ref{cor:3.6}]\label{thm:1.3}
Let $\Gamma$ be a Polish group, and let $\mu$ be a Borel probability measure on $\Gamma$. Let $S\subseteq\Gamma$ be a countably infinite subset. Then the Schreier graph $\Sch(\Gamma,S,\Gamma)$ admits a $\mu$-measurable domatic $\aleph_0$-partition. 
\end{theorem}

We now explain how these results above relate to Zelinka's theorem \citep{zelinka} on finite hypercube graphs introduced earlier. Let $Q_{\N}$ be the loop-free simple undirected graph on the vertex set $V=\{0,1\}^\N$ of all infinite binary sequences, such that $(v,w)\in V^2$ is an edge of $Q_{\N}$ if and only if the two infinite binary sequences $v,w\in\{0,1\}^{\N}$ differ exactly in one place. Thus $Q_{\N}$ is the $\aleph_0$-dimensional version of the finite hypercube graphs $Q_n$. 

The graph $Q_{\N}$ is then isomorphic to a Schreier graph $\Sch(\Gamma,S,\Gamma)$, where $\Gamma=(\Z/2\Z)^{\N}$ is a zero-dimensional compact Polish group, and $S$ is the set of all sequences $s\in\Gamma$ which contains a $1\in\Z/2\Z$ in exactly one place and $0\in\Z/2\Z$ elsewhere. Since the topological closure $\overline{S}\subseteq\Gamma$ is countable, Theorem \ref{thm:1.1} implies that the graph $Q_{\N}\cong\Sch(\Gamma,S,\Gamma)$ does not admit any Borel or Baire measurable domatic $\aleph_0$-partition. Since $S\subseteq \Gamma$ is countably infinite, Theorem \ref{thm:1.3} implies that $Q_{\N}\cong\Sch(\Gamma,S,\Gamma)$ does admit measure-theoretic domatic $\aleph_0$-partitions. 

The main results of Section \ref{sec:2}, Theorems \ref{thm:1.1} and \ref{thm:1.2}, are proved using a main lemma on the existence of domatic finite partitions with open parts, which extends a theorem in Alon--Spencer \citep[Theorem 5.2.2]{alon-spencer}. 

\begin{lemma}[Theorem \ref{thm:2.12}]\label{lem:1.4}
Let $\Gamma$ be a locally compact Polish group with a two-sided invariant metric (eg.\ when $\Gamma$ is compact or abelian \citep[\S{}2.1]{gao}) and finite topological dimension. For every $k,n\in\N$, there exists some $N=N(k,n)\in\N$, such that for any sets $F_0,\ldots,F_{n-1}\subseteq\Gamma$ with $\tabs{F_0}=\ldots=\tabs{F_{n-1}}=N$, there exists a sequence of pairwise-disjoint open subsets $D_0,\ldots,D_{k-1}$ of $\Gamma$, for which every right translate $F_i\cdot \gamma$ of every $F_i$ intersects every set $D_j$. In particular, the sequence $\langle D_j:j<k\rangle$ is a domatic $k$-partition with open parts for each of the graphs $\Sch(\Gamma,F_0,\Gamma),\ldots,\Sch(\Gamma,F_{n-1},\Gamma)$. 
\end{lemma}

In Section \ref{sec:2.6}, we find an application of our analysis to the theory of sum sets, as we give an extension of a theorem by Erd\H{o}s--Kunen--Mauldin \citep[Theorem 1]{ekm}:

\begin{theorem}[Corollary \ref{cor:2.29}]\label{thm:1.5}
Let $1\le n\in\N$, and let $P\subseteq\R^n$ be a nonempty closed perfect subset of $\R^n$. Then there exists a family $\langle C_i:i<2^{\aleph_0}\rangle$ of $2^{\aleph_0}$ pairwise-disjoint closed subsets of $\R^n$, such that $P+C_i=\R^n$ and $C_i+C_j=\R^n$ for all $i,j<2^{\aleph_0}$. 
\end{theorem}

In Sections \ref{sec:3} and \ref{sec:4}, we list out other results concerning domatic partitions for Borel graphs in general. Section \ref{sec:3} mainly concerns the existence of domatic $\aleph_0$-partitions on $\aleph_0$-regular Borel graphs, and Section \ref{sec:4} mainly concerns the existence or nonexistence of domatic finite partitions on locally countable Borel graphs. Notable results include:

\begin{theorem}[Theorem \ref{thm:3.5}]\label{thm:1.6}
Let $G$ be a an out-degree $\aleph_0$-regular Borel graph with countable in-degrees on a Borel probability space $(X,\mu)$ of vertices. Then $G$ admits a $\mu$-measurable domatic $\aleph_0$-partition. 
\end{theorem}

\begin{theorem}[Theorem \ref{thm:4.3}]\label{thm:1.7}
There exists a fully looped undirected $\aleph_0$-regular acyclic Borel graph $G$ on a Polish space $(X,\tau)$ of vertices, without $\tau$-Baire measurable domatic $3$-partitions. 
\end{theorem}

The next two results concern an edge-coloring version of domatic partitions; see Definition \ref{def:3.7}. 

\begin{theorem}[Theorem \ref{thm:3.8}]\label{thm:1.8}
Let $G$ be a loop-free simple undirected $\aleph_0$-regular Borel graph on a Borel space $X$ of vertices. If $\mu$ is any Borel probability measure on $X$, then there is an $E_G$-invariant $\mu$-conull Borel set $C_\mu\subseteq X$ and a symmetric Borel function $f_\mu:G\restrict{}C_\mu\to \N$ such that $f_\mu$ is domatic everywhere in $C_\mu$. Similarly, if $\tau$ is any Polish topology on $X$, then there is an $E_G$-invariant $\tau$-comeager Borel set $C_\tau\subseteq X$ and a symmetric Borel function $f_\tau:G\restrict{} C_\tau\to \N$ such that $f_\tau$ is domatic everywhere in $C_\tau$. 
\end{theorem}

\begin{theorem}[Weilacher, Theorem \ref{thm:4.5}]\label{thm:1.9}
There exists a loop-free simple undirected $\aleph_0$-regular acyclic Borel graph $G\subseteq X^2$ on a Borel space $X$ of vertices, without symmetric Borel domatic edge-$2$-partitions. Moreover, $G$ is Borel bipartite without Borel sinkless orientations. 
\end{theorem} 

Finally in Section \ref{sec:4.4}, we discuss the question of how much can be said of the existence or nonexistence of domatic finite partitions on locally finite Borel graphs. We also leave an open question in that section: Let $\Gamma$ be a countably infinite group and let $S\subseteq \Gamma$ be an arbitrary countably infinite generating set. Is it necessarily true that for every free Borel $\Gamma$-space $X$, the out-degree $\aleph_0$-regular Borel graph $\Sch(\Gamma,S,X)=\{(x,y)\in X^2:\exists s\in S\left(s\cdot x=y\right)\}$ on $X$ admits Borel domatic $k$-partitions for every finite $k\in\N$?

\subsection{Notation}\label{sec:1.1}

Let $G$ be a directed graph with possible loops on a vertex set $X$, represented as a binary relation $G\subseteq X^2$. We define the \emph{out neighborhood} (or simply just \emph{neighborhood}) of a vertex $x\in X$ to be the set $N_G(x)=\{y\in X:(x,y)\in G\}$. We define the \emph{out-degree} (or just \emph{degree}) of $x\in X$ to be the cardinality $\tabs{N_G(x)}$, and we say that $G$ is \emph{$\kappa$-regular} if every vertex has out-degree $\kappa$. 

We write $E_G\subseteq X^2$ for the connectedness equivalence relation of a graph $G$, and say a set $A\subseteq X$ is \emph{$E_G$-invariant} if $A$ is closed under $E_G$-equivalence, or equivalently if $A$ is a union of $G$-connected components. A \emph{coloring} is a function $f:X\to Y$ such that $f(x)\ne f(y)$ for all edges $(x,y)\in G$ where $x\ne y$. When talking about partial functions $f:X\partialto Y$ whose domain is a subset of the vertex set $X$, we will frequently call the codomain $Y$ the set of \emph{colors} of $f$. We say a set $I\subseteq X$ is \emph{independent} if for every edge $(x,y)\in G$ with $x\ne y$, not both $x$ and $y$ belong to $I$. 

The diagonal set $\Delta_X=\{(x,x)\in X^2: x\in X\}$ is the set of \emph{loops} on the vertex set $X$. We say the graph $G$ is \emph{loop-free} if $G\cap \Delta_X=\varnothing$, and we say that $G$ is \emph{fully looped} if $\Delta_X\subseteq G$. Thus the graphs $G\ssm \Delta_X$ and $G\cup \Delta_X$ are the loop-free and fully looped versions of $G$ respectively. Note that if $G$ is fully looped, then $x\in N_G(x)$ for every $x\in X$. We say that an undirected graph with loops $G$ is \emph{acyclic} if its loop-free version $G\ssm \Delta_X$ is acyclic. 

A set $D\subseteq X$ of vertices is \emph{dominating} for $G$ if it intersects every neighborhood set, meaning $D\cap N_G(x)\ne\varnothing$ for all $x\in X$. If $\kappa$ is a cardinal number, a \emph{domatic $\kappa$-partition} for $G$ is a sequence of $\kappa$ pairwise-disjoint dominating sets. A partial function $f:X\partialto \kappa$ is \emph{domatic} at a vertex $x$ if $f[N_G(x)]=\kappa$, and $f$ is \emph{domatic} if it's domatic everywhere. From a domatic partial function $f:X\partialto\kappa$ with $\kappa$ colors, one can produce a domatic $\kappa$-partition $\langle f^{-1}[\{i\}]\subseteq X:i<\kappa\rangle$, and thus a domatic $\kappa$-partition is equivalent to a domatic partial function with $\kappa$ colors. If $1\le \kappa\le \aleph_0$ and $\mathcal{F}$ is a $\sigma$-algebra on $X$, then $G$ admits an $\mathcal{F}$-measurable domatic partial $\kappa$-partition if and only if $G$ admits an $\mathcal{F}$-measurable domatic total $\kappa$-partition, since one can paint all uncolored vertices in $X$ with a fixed junk color in $\kappa$. 

We note here the monotonicity of domaticity. If $G\subseteq H$ are graphs on a same vertex set $X$, then $G$ admitting a domatic $\kappa$-partition implies that $H$ admits the same domatic $\kappa$-partition. Similarly if $\kappa\le \lambda$ are cardinals, then a graph $G$ admitting a domatic $\lambda$-partition implies that it admits a domatic $\kappa$-partition. Intuitively, the more edges a graph has, the smaller its dominating sets become, and the easier it is to pack more dominating sets into its set of vertices. 

An important class of Borel graphs in descriptive graph combinatorics is the Schreier graphs. When $\Gamma$ is a group acting on a set $X$ and $S\subseteq\Gamma$ is a (possibly not generating) set, we associate the directed \emph{Schreier graph} $\Sch(\Gamma,S,X)=\{(x,y)\in X^2:\exists s\in S\left(s\cdot x=y\right)\}$ over the vertex set $X$. When $\Gamma$ is Polish, $X$ is Borel, and the action of $\Gamma$ on $X$ is Borel, the graph $\Sch(\Gamma,S,X)$ is also Borel when $S$ is countable, or when $S$ is Borel and $\Gamma$ acts freely on $X$. Note that the out-neighborhood sets of the Schreier graph are given by $N_{\Sch(\Gamma,S,X)}(x)=S\cdot x$ for all $x\in X$. 

The symbol $\omega$ is used to mean the least infinite ordinal, and it is equal to the set $\omega=\N=\{0,1,2,\ldots\}$ of nonnegative integers. The symbol $\infty$ means the positive infinity in the extended real numbers. Thus the expression ``$n<\omega$'' means $n$ is a nonnegative integer, and the expression ``$r<\infty$'' means $r$ is a finite real number. The phrase ``perfect set'' is used by default to mean a closed perfect set inside a Polish space, which distinction will not matter except when it's explicitly disambiguated. 

\subsection*{Acknowledgements}

We thank Clinton Conley for his constant help and many insightful conversations in the writing of this paper. We thank Felix Weilacher for proving Theorem \ref{thm:4.5} in private communications, a cleaner presentation of Theorem \ref{thm:3.5}, and many other insightful feedback. We thank an anonymous referee for helpful feedback and a cleaner proof of Theorem \ref{thm:2.1}. We thank Alexander Kechris for helpful comments on the presentation of this paper. 

\section{Proof of Theorem \ref{thm:1.1}}\label{sec:2}

The main goal of the following few sections is to prove Corollaries \ref{cor:2.18} and \ref{cor:2.19}, which concern Borel and Baire measurable domatic $\aleph_0$-partitions of Schreier graphs defined from finite-dimensional compact Polish group actions. 

More precisely, we fix an infinite compact Polish group $\Gamma$ with finite topological dimension and an arbitrary subset $S\subseteq \Gamma$. We fix the continuous action of $\Gamma$ on itself via left multiplication, and we recall that the Schreier graph $G=\Sch(\Gamma,S,\Gamma)$ on $\Gamma$ is defined by $G=\{(\gamma,s\cdot\gamma):\gamma\in\Gamma,s\in S\}$. 

In Section \ref{sec:2.1} we prove if $\overline{S}\subseteq\Gamma$ is countable compact then $\Sch(\Gamma,S,\Gamma)$ does not admit Baire measurable domatic $\aleph_0$-partitions. 

Sections \ref{sec:2.2} to \ref{sec:2.4} are fully devoted to proving one technical black box Theorem \ref{thm:2.12}. In Section \ref{sec:2.2} we give a characterization of the topological dimension of $\Gamma$ using the Gleason--Yamabe theorem. Section \ref{sec:2.3} uses the finite dimension of $\Gamma$ to construct a packing of $\Gamma$ with open cells, and Section \ref{sec:2.4} uses finite dimension again to show that a random finite coloring of these open cells (in the sense of the Lov\'{a}sz local lemma) gives domatic finite partitions on $\Gamma$ with open parts. In Section \ref{sec:2.5}, we use the compactness of $\Gamma$ to show if $\overline{S}\subseteq\Gamma$ is uncountable then domatic finite partitions on $\Gamma$ with open parts can always be lifted to domatic $\aleph_0$-partitions. 

Finally in Section \ref{sec:2.6}, we present an application of our methods to prove a result Corollary \ref{cor:2.29} about sum sets in $\R^n$. 

\subsection{Countable compactness implies anti-domaticity}\label{sec:2.1}

\begin{theorem}\label{thm:2.1}
Let a Polish group $\Gamma$ continuously act on a Polish space $X$. Let $S\subseteq\Gamma$ be a countable compact set, with its Schreier graph $G\coloneqq\Sch(\Gamma,S,X)$ on $X$. For any Baire measurable function $f:X\to\omega$, there is a comeager set of $x\in X$ for which $f[N_{G}(x)]$ is finite. In particular, $f$ is not domatic at any such vertex $x$. 
\end{theorem}

\begin{proof}
By Kechris \citep[Theorem 8.38]{kechris}, we can fix some comeager $G_{\delta}$ set $A\subseteq X$ such that $f\restrict{}A$ is continuous. Since every $\gamma\in \Gamma$ acts on $X$ by a homeomorphism, and hence preserves comeager-ness of $A$, the intersection $A'=\bigcap_{\gamma\in \langle S\rangle}\gamma\cdot A\subseteq A$ over the countable subgroup $\langle S\rangle$ generated by $S$ is also comeager $G_\delta$. Then $f\restrict{}A'$ is also continuous, and $A'$ is $E_G$-invariant. 

For each $x\in A'$, since $A'$ is $E_{\Sch(\Gamma,S,X)}$-invariant, we have $S\cdot x\subseteq A'$ and so $f\restrict{}A'$ is continuous over $S\cdot x$. The function $g:S\to \omega$ defined by $g(s)=f(s\cdot x)$ is continuous as it is a composition of the continuous functions $s\mapsto s\cdot x$ and $f\restrict{}A'$. Since $S$ is compact, its continuous image $g[S]\subseteq \omega$ must be finite. Thus we have shown that for all $x$ inside the comeager set $A'$, the set $f[N_G(x)]=f[S\cdot x]=g[S]\subseteq \omega$ is finite, as desired. 
\end{proof}

Since the nonexistence of domatic partitions can be passed to subgraphs, the same result will hold if the set $S$ is only assumed to have countable compact topological closure in $\Gamma$. 

\subsection{The dimension of a locally compact Polish group}\label{sec:2.2}

\begin{theorem}[Gleason--Yamabe, see {\citep{tao}}]\label{thm:2.2}
Let $G$ be a locally compact group. Then, for any open neighbourhood $U$ of the identity, there exists an open subgroup $G'$ of $G$ and a compact normal subgroup $K$ of $G'$ in $U$ such that $G'/K$ is isomorphic to a Lie group.
\end{theorem}

\begin{corollary}\label{cor:2.3}
Let $\Gamma$ be a locally compact Polish group. Then $\Gamma$ is an inverse limit of an inverse system $\langle \Gamma_i:i<\omega\rangle$ of Lie groups and continuous surjective homomorphisms:
\[\Gamma=\varprojlim\nolimits_i \Gamma_i\to\;\cdots\;\to \Gamma_2\to\Gamma_1\to\Gamma_0\]
\end{corollary}

\begin{proof}
Let $\{U_i:i<\omega\}$ be an open neighborhood basis at $1_\Gamma\in\Gamma$. Fix a sequence $K_0,K_1,\ldots$ of compact normal subgroups of $G$, such that $K_i\subseteq U_i$ and $\Gamma/K_i$ is a Lie group for all $i<\omega$. The group $\Gamma_i\coloneqq \Gamma/(K_0\cap\ldots\cap K_i)$ is still a Lie group for $i<\omega$, since it embeds as a closed subgroup into the Lie group $(\Gamma/K_0)\times\ldots\times(\Gamma/K_i)$. The inverse system $\langle \Gamma_i\rangle$ is as desired. 
\end{proof}

\begin{definition}\label{def:2.4}
Let $\Gamma$ be a locally compact Polish group, and fix any inverse system $\langle \Gamma_i:i<\omega\rangle$ associated with $\Gamma$ as above. We define the \emph{dimension} of $\Gamma$ to be the supremum of the dimensions of the real manifolds $\Gamma_i$, written as $\dim(\Gamma)\coloneqq\sup\{\dim_\R(\Gamma_i):i<\omega\}<\omega+1$. We see that $\dim(\Gamma)$ can take any value from $\{0,1,2,\ldots,\aleph_0\}$. 
\end{definition}

The way we'll use the dimension as we defined is via the next two lemmas. 

\begin{lemma}\label{lem:2.5}
Let $\Gamma$ be a locally compact Polish group of dimension $\dim(\Gamma)$. Then $\dim(\Gamma)$ is equal to the small inductive dimension, the large inductive dimension, and the Lebesgue covering dimension of $\Gamma$ as a topological space. In particular, $\dim(\Gamma)$ is well-defined. 
\end{lemma}

\begin{proof}
By the Kat\v{e}tov--Morita theorem \citep[Theorem 7.3.3]{engelking}, the three topological dimensions of a separable metrizable space agree. We will finish by showing that $\dim(\Gamma)=d(\Gamma)$, where $d(\Gamma)$ denotes the Lebesgue covering dimension of $\Gamma$. 

First we'll show for all finite $n<\omega$, if $\dim(\Gamma)\le n$ then $d(\Gamma)\le n$. If $\dim(\Gamma)\le n$, the each $\Gamma_i$ is a manifold of dimension $\le n$, and so each $\Gamma_i$ has Lebesgue covering dimension $d(\Gamma_i)\le n$. By Nagami's theorem \citep[Exercise 7.3.I]{engelking}, their inverse limit $\Gamma=\varprojlim_i\Gamma_i$ also has $d(\Gamma)\le n$. 

Next we'll show for all finite $n<\omega$, if $\dim(\Gamma)\ge n$ then $d(\Gamma)\ge n$, which completes the proof. If $\dim(\Gamma)\ge n$, then there is some $\Gamma_i$ whose manifold dimension $d_i$ satisfies $d_i\ge n$. For each $j<\omega$, let $\exp_j:\mathfrak{g}_j\to \Gamma_j$ be the exponential map associated with $\Gamma_j$. One can check that taking the inverse limit as $j\to\infty$ gives a well-defined exponential map $\exp:\mathfrak{g}\to\Gamma$, and that the projection $\pi:\Gamma\to\Gamma_i$ induces a projection $\pi_*:\mathfrak{g}\to\mathfrak{g}_i$. Note that $\mathfrak{g}$ and $\mathfrak{g}_i$ are linearly isomorphic to $\mathfrak{g}\cong\R^{\dim(\Gamma)}$ and $\mathfrak{g_i}\cong\R^{d_i}$ respectively. 

Let $U\subseteq \mathfrak{g}_i$ be an open neighborhood of $0\in\mathfrak{g}_i$ such that the exponential map $\exp_i:\mathfrak{g}_i\to\Gamma_i$ induces a homeomorphism $U\cong\exp_i(U)$. We lift $U$ linearly under the projection $\pi_*:\mathfrak{g}\to\mathfrak{g_i}$ to a $d_i$-manifold $0\in V\subseteq\mathfrak{g}$ such that $\pi_*:V\cong U$ is a homeomorphism. The composite map $\pi\circ \exp={\exp_i}\circ {\pi_*}:\mathfrak{g}\to\Gamma_i$ thus induces a homeomorphism $V\cong\exp_i(U)$, and so its factor map $\pi$ induces the homeomorphism $\exp(V)\cong\exp_i(U)$. 

Since $U\subseteq\mathfrak{g}_i\cong\R^{d_i}$ has closed subsets of Lebesgue covering dimension $d_i$, its homeomorphic copy $\exp(V)\subseteq\Gamma$ also has closed subsets of Lebesgue covering dimension $d_i$. Since the Lebesgue covering dimension is hereditary to closed subsets \citep[Theorem 7.1.8]{engelking}, we see that $\Gamma\supseteq\exp(V)$ has Lebesgue covering dimension $d(\Gamma)\ge d_i\ge n$ as desired. {}
\end{proof}

\begin{lemma}\label{lem:2.6}
Let $\Gamma$ be a locally compact Polish group with $\dim(\Gamma)<\infty$. Then there is a finite positive constant $M_\Gamma<\infty$ and a neighborhood basis of the identity $1_\Gamma$ consisting of open sets $U$ such that every family of left translates $\gamma\cdot U$ of $U$ that are pairwise-disjoint subsets of $(UU^{-1})^4U$ has at most $M_\Gamma$ many members. 
\end{lemma}

\begin{proof}
Since $\Gamma$ is an inverse limit of Lie groups of dimension at most $\dim(\Gamma)<\infty$, it suffices to show that every $d$-dimensional Lie group $G$ has such a basis with $M_G=10^d$. 

Since $G$ is a Lie group, it admits left-invariant Riemannian metrics, and so it also admits a complete left-invariant metric $d_G$ such that for any $\varepsilon>0$ there exists some open neighborhood $1_G\in\Omega\subseteq G$ of the identity and a $(1+\varepsilon)$-bi-Lipschitz homeomorphism $f:D\cong \Omega$ from some open Euclidean domain $D\subseteq\R^d$ to $\Omega$. 

Then, for all sufficiently small $r>0$, the open ball $U=B_{d_G}(1_G,r)\subseteq G$ is such that $U=U^{-1}$, there is $U^9\subseteq\Omega$ and $f^{-1}[U^9]\subseteq D$ is contained in a ball of radius $9(1+\varepsilon)r$, and every left translate $\gamma\cdot U\subseteq U^9$ is a radius-$r$ ball $\gamma\cdot U=B_{d_G}(\gamma,r)$ such that $f^{-1}[\gamma U]\subseteq f^{-1}[U^9]$ contains some ball of radius $(1-\varepsilon)r$. For sufficiently small $\varepsilon>0$, a Euclidean ball of radius $9(1+\varepsilon)r$ has volume at most $10^d$ times that of radius $(1-\varepsilon)r$, which means at most $10^d$ many left translates of $U$ can be packed into $U^9=(UU^{-1})^4U$. {}
\end{proof}

\subsection{An open-cells packing}\label{sec:2.3}

In this section, we prove Theorem \ref{thm:2.9}. 

\begin{lemma}\label{lem:2.7}
Assume $X$ is a metrizable space, $A\subseteq X$ is closed, $B\subseteq X$ is open such that $A\subseteq B$, $C\subseteq X$ is closed, and $U\subseteq C$ is relatively open in $C$ such that $A\cap C\subseteq U\subseteq \overline{U}\subseteq B\cap C$. Then there is a set $V\subseteq X$ open in $X$ such that $A\subseteq V\subseteq \overline{V}\subseteq B$, $V\cap C=U$, and $\overline{V}\cap C=\overline{U}$. 
\end{lemma}

\begin{proof}
Fix a metric $d$ on $X$, and let $W$ be an open set in $X$ such that $A\cup\overline{U}\subseteq W\subseteq\overline{W}\subseteq B$. We define $V=\{x\in W:d(x,A\cup\overline{U})<d(x,C\ssm U)\}$ as desired. {}
\end{proof}

In the following proofs, we use the phrase ``topological dimension'' or ``$\dim(X)$'' of a Polish space $X$ to mean any of its small or large inductive dimension, or its Lebesgue covering dimension. These dimensions are all equal by the Kat\v{e}tov--Morita theorem \citep[Theorem 7.3.3]{engelking}.

\begin{lemma}\label{lem:2.8}
Let $X$ be a Polish space. Assume $r<\omega$, and $M_0,\ldots,M_{r-1}$ are nonempty closed subsets of $X$ of finite topological dimensions $0\le \dim(M_i)<\infty$. If $A\subseteq X$ is closed and $B\subseteq X$ is open such that $A\subseteq B$, then there exists a set $U\subseteq X$ open in $X$ such that $A\subseteq U\subseteq\overline{U}\subseteq B$, and for all $i<r$, we have $\dim(\partial U\cap M_i)<\dim(M_i)$. 
\end{lemma}

\begin{proof}
We proceed by induction on $r<\omega$, noting that there is nothing to prove when $r=0$. Assume we have proved the lemma for a fixed $r<\omega$, and we will prove next the case with $r+1$ given closed sets $M_0,\ldots,M_r\subseteq X$. We may assume $\dim(M_0)\le\ldots\le \dim(M_r)$. 

By the inductive hypothesis, let $U$ be an open set such that $A\subseteq U\subseteq\overline{U}\subseteq B$, and such that $\dim(\partial U\cap M_i)<\dim(M_i)$ for $0\le i\le r-1$. Let $M'=M_0\cup\ldots\cup M_{r-1}$. We will next work in the subspace $Y=(M'\cup M_r)\ssm(\partial U\cap M')$. Note that by the countable sum theorem \citep[Theorem 7.2.1]{engelking} and by that the topological dimension is hereditary to subspaces \citep[Theorem 7.1.1]{engelking}, we have $\dim(Y)\le\dim(M'\cup M_r)\le\dim(M_r)$ and $\dim(\partial U\cap M')<\dim(M_r)$.

Let $Z=U\cap M'$, and note that $Z$ is relatively closed in $Y$ since $Z=\overline{U}\cap M'\cap Y$. Then we have $(A\cap Y)\cup Z\subseteq U\cap Y$, where $(A\cap Y)\cup Z$ is relatively closed in $Y$, and $U\cap Y$ is relatively open in $Y$. By definition of the large inductive dimension \citep[\S{}7.1]{engelking}, there exists a set $V\subseteq Y$ relatively open in $Y$ such that $(A\cap Y)\cup Z\subseteq V\subseteq U\cap Y$, and $\dim(\partial^Y V)<\dim(Y)\le\dim(M_r)$ where $\partial^Y V=(\overline{V}\ssm V)\cap Y$ is the relative boundary of $V$ in $Y$. 

Since $\partial U\cap M'$ is closed, $Y$ is relatively open in $M'\cup M_r$. Let $C=M'\cup M_r$. Since $V$ is relatively open in $Y$, we see that $V$ is also relatively open in $C$, and we see that the relative boundary of $V$ in $C$ is $\partial^{C}V\subseteq(\partial^Y V)\cup (\partial U\cap M')$. We have the following by the countable sum theorem \citep[Theorem 7.2.1]{engelking}, since $\partial^Y V$ is $F_\sigma$ and $\partial U\cap M'$ is closed:
\[\dim(\partial^C V\cap M_r)\le \dim(\partial^C V)\le \max\{\dim(\partial^Y V),\dim(\partial U\cap M')\}<\dim(M_r)\]

We also have the following:
\[A\cap C=A\cap U\cap C=A\cap Y\subseteq V\subseteq\overline{V}\subseteq\overline{U\cap Y}\subseteq\overline{U}\cap C\subseteq B\cap C\]
Then, since $U\cap M'=Z\subseteq V\cap M'$ and $\overline{V}\cap M'\subseteq \overline{U}\cap C\cap M'=\overline{U}\cap M'$, we have the inclusion $\partial^{C}V\cap M'\subseteq\partial U\cap M'$, which means by the inductive hypothesis that for all $0\le i\le r-1$, we have $\dim(\partial^C V\cap M_i)\le \dim(\partial U\cap M_i)<\dim(M_i)$. 

Finally, by Lemma \ref{lem:2.7}, there is an open set $W\subseteq X$ extending $V\subseteq C$ such that $A\subseteq W\subseteq\overline{W}\subseteq B$ and $\partial W\cap C=\partial^C V$. This $W$ completes the inductive step for the case $r+1$ as desired. {}
\end{proof}

\begin{theorem}\label{thm:2.9}
Let $\Gamma$ be a Polish group with finite topological dimension $d<\infty$ as a Polish space. Let $F\subseteq\Gamma$ be a finite set of size $\tabs{F}=n\ge d$, and let $\{U_i:i<\omega\}$ be an open cover of $\Gamma$ such that every $U_i$ satisfies ${(FF^{-1}\ssm\{1_\Gamma\})\cdot U_i}\cap {U_i}=\varnothing$ (eg.\ when $\diam(U_i)$ is sufficiently small). Then there is a family $\mathcal{R}=\{R_i:i<\omega\}$ of pairwise-disjoint open subsets of $\Gamma$, such that $R_i\subseteq U_i$ for all $i<\omega$, and moreover for every $\gamma\in \Gamma$, the right translate $F\cdot \gamma$ of $F$ intersects at least $n-d$ many distinct members of $\mathcal{R}$, i.e.\ $\tabs{\{i<\omega:F\cdot\gamma\cap R_i\ne\varnothing\}}\ge n-d$. 
\end{theorem}

\begin{proof}
Since every open set is $F_\sigma$, there is a countable closed covering $\{A_i:i<\omega\}$ refining the open covering $\{U_i:i<\omega\}$. We can fix a function $a:\omega\to \omega$ such that $A_i\subseteq U_{a(i)}$ for all $i<\omega$. 

The first step of this proof is to construct a sequence $\langle V_i:i<\omega\rangle$ of open sets $V_i$ inductively on $i<\omega$, such that in constructing each set $V_i$ we ensure that $A_i\subseteq V_i\subseteq\overline{V_i}\subseteq U_{a(i)}$, and that for every $f\in F$ and $S\subseteq\{0,\ldots,i-1\}$, if the closed set $M_{f,S}=\bigcap_{j\in S}fF^{-1}\cdot\partial V_j$ is nonempty, then $\dim(\partial V_i\cap M_{f,S})<\dim(M_{f,S})$. (Note that when $S=\varnothing$, the last condition just says $\dim(\partial V_i)<\dim(\Gamma)=d$.) Since for each $i<\omega$, the number of such pairs $(f,S)$ is finite, the construction of $\langle V_i:i<\omega\rangle$ follows from an application of Lemma \ref{lem:2.8}. 

We now prove by induction on $0\le t\le d+1$ that for every sequence $s_0<\ldots<s_{t-1}<\omega$, we have $\dim(F^{-1}\cdot\partial V_{s_0}\cap\ldots\cap F^{-1}\cdot\partial V_{s_{t-1}})\le d-t$. When $t=0$, this says $\dim(\Gamma)\le d$. In the inductive step going from case $t$ to case $t+1$, we have the following:
\begin{align*}
\dim\left(\bigcap_{{i\le t}}F^{-1}\partial V_{s_i}\right)&=\dim\left(\bigcup_{{f\in F}}\left(\bigcap_{{i\le t-1}}F^{-1}\partial V_{s_i}\cap f^{-1}\partial V_{s_t}\right)\right)\\
&\le \max_{f\in F}\dim\left(\bigcap_{{i\le t-1}}F^{-1}\partial V_{s_i}\cap f^{-1}\partial V_{s_t}\right)\\
&=\max_{f\in F}\dim\left(\bigcap_{i\le t-1}fF^{-1}\partial V_{s_i}\cap \partial V_{s_t}\right)\\
&\le\max_{f\in F}\max\left\{\dim\left(\bigcap_{i\le t-1}fF^{-1}\partial V_{s_i}\right)-1,\;-1\right\}\\
&=\max\left\{\dim\left(\bigcap_{i\le t-1}F^{-1}\partial V_{s_i}\right)-1,\;-1\right\}\\
&\le d-t-1
\end{align*}
Here we used the sum theorem for dimension \citep[Theorem 7.2.1]{engelking}, that $F$ act on $\Gamma$ by homeomorphisms, the construction of $V_{s_t}$, and the inductive hypothesis. Thus in particular when $t=d+1$, we get that for every sequence $s_0<\ldots<s_d<\omega$, we have $F^{-1}\partial V_{s_0}\cap\ldots\cap F^{-1}\partial V_{s_d}=\varnothing$. 

Next, we claim that any right translate $F\cdot\gamma$ of $F$ can intersect $\bigcup_{i<\omega}\partial V_i$ at most $d$ times. If not, then there exists a $\gamma\in \Gamma$ and distinct elements $f_0,\ldots,f_d\in F$ such that $f_i\cdot\gamma\in \bigcup_{s<\omega}\partial V_s$ for all $i\le d$. For each $i$, we can pick some $s_i<\omega$ such that $f_i\cdot \gamma\in \partial V_{s_i}\subseteq U_{a(s_i)}$. Then for $i\ne j$, since $f_jf_i^{-1}U_{a(s_i)}\cap U_{a(s_i)}=\varnothing$ by assumption and $f_j\cdot\gamma\in f_jf_i^{-1}U_{a(s_i)}\cap U_{a(s_j)}\ne\varnothing$, we have $a(s_i)\ne a(s_j)$, and hence $s_0,\ldots,s_d$ are all distinct. Then $\gamma\in F^{-1}\partial V_{s_0}\cap\ldots\cap F^{-1}\partial V_{s_d}\ne\varnothing$, which contradicts our earlier arguments. 

Let $W_i=V_i\ssm\bigcup_{j<i}\overline{V_j}$ for $i<\omega$. Since $\bigcup_{i<\omega}V_i\supseteq \bigcup_{i<\omega}A_i=\Gamma$, we see that $\bigcup_{i<\omega}W_i\cup\bigcup_{i<\omega}\partial V_i=\Gamma$, while $\{W_i:i<\omega\}$ is a family of pairwise-disjoint open sets. For $i<\omega$, if we let $R_i=\bigcup_{a(j)=i}W_j$, then $\mathcal{R}=\{R_i:i<\omega\}$ is a family of pairwise-disjoint open sets such that $R_i\subseteq U_i$ for all $i$. The previous argument shows that every right translate $F\cdot\gamma$ of $F$ intersects the set $\bigcup_{i<\omega}R_i\supseteq\Gamma\ssm\bigcup_{i<\omega}\partial V_i$ at least $n-d$ times, and a similar argument as before using our assumptions on $U_i$ shows that in each of these $n-d$ times, $F\cdot\gamma$ must intersect a distinct $R_i\in\mathcal{R}$. The family $\mathcal{R}=\{R_i:i<\omega\}$ is as desired. {}
\end{proof}

\subsection{An open cover of locally bounded growth}\label{sec:2.4}

\begin{lemma}\label{lem:2.10}
Assume $\Gamma$ is a Polish group with a two-sided invariant metric (eg.\ when $\Gamma$ is compact or abelian \citep[\S{}2.1]{gao}) and, for some absolute constant $M<\infty$, a neighborhood basis of the identity $1_\Gamma$ consisting of open sets $U$, at most $M$ of whose left translates $\gamma\cdot U$ can be packed into $(UU^{-1})^4U$. Assume $F\subseteq\Gamma$ is a finite set. Then there is an open cover $\mathcal{B}$ of $\Gamma$ such that for every $V\in \mathcal{B}$, $(FF^{-1}\ssm\{1_\Gamma\})\cdot V\cap V=\varnothing$, every left translate $\gamma\cdot V$ intersects at most $M$ members of $\mathcal{B}$, there are no distinct $f,g\in F$ for which both $fV\cap A\ne\varnothing$ and $gV\cap A\ne\varnothing$ for some $A\in\mathcal{B}$, and the number of $W\in\mathcal{B}$ such that both $FV\cap A\ne \varnothing$ and $FW\cap A\ne\varnothing$ for some $A\in\mathcal{B}$ is at most $M\cdot\tabs{F}^2$. 
\end{lemma}

\begin{proof}
Since $\Gamma$ admits a two-sided invariant metric, $1_\Gamma$ is not a limit point of the conjugacy classes of $FF^{-1}\ssm\{1_\Gamma\}$. We can fix an open neighborhood $U$ of the identity $1_\Gamma$ such that $(UU^{-1})^4$ is disjoint from the conjugacy classes of $FF^{-1}\ssm\{1_\Gamma\}$, and at most $M$ left translates $\gamma\cdot U$ can be packed into $(UU^{-1})^4U$. 

Let $S\subseteq\Gamma$ be a maximal set such that the family $\{s\cdot U:s\in S\}$ is pairwise-disjoint. The family $\mathcal{B}=\{s\cdot UU^{-1}:s\in S\}$ is an open cover of $\Gamma$, since for all $x\in\Gamma$, if $s\in S$ is such that $sU\cap xU\ne\varnothing$ then $x\in sUU^{-1}\in\mathcal{B}$. We'll next check that $\mathcal{B}$ has our desired properties. 

Let $V=sUU^{-1}\in \mathcal{B}$ and $g\in FF^{-1}\ssm\{1_\Gamma\}$. Then $g\cdot sUU^{-1}\cap sUU^{-1}=\varnothing$ follows from our assumption that $s^{-1}gs\notin (UU^{-1})^2\subseteq (UU^{-1})^4$. This means $(FF^{-1}\ssm\{1_\Gamma\})\cdot V\cap V=\varnothing$. 

Fix $\gamma\in \Gamma$, and let $S_\gamma\subseteq S$ be the set of all $s\in S$ for which $\gamma UU^{-1}\cap sUU^{-1}\ne\varnothing$. Then we see that $\gamma^{-1}s\in(UU^{-1})^2$ for all $s\in S_\gamma$, and since $S_\gamma\subseteq S$, the family $\{\gamma^{-1}sU:s\in S_\gamma\}$ is a pairwise-disjoint family of left translates of $U$ which are subsets of $(UU^{-1})^2U\subseteq(UU^{-1})^4U$, so by assumption we get $\tabs{S_\gamma}\le M$. Since every $V\in\mathcal{B}$ is a left translate of $UU^{-1}$, we get that every left translate of $V$ also intersects at most $M$ members $sUU^{-1}\in\mathcal{B}$. 

Let $V=xUU^{-1}\in\mathcal{B}$, and let $f,g\in F$ be distinct such that there is some $A=sUU^{-1}\in \mathcal{B}$ for which both $fV\cap A\ne\varnothing$ and $gV\cap A\ne\varnothing$. From $fxUU^{-1}\cap sUU^{-1}\ne\varnothing$ we get $x^{-1}f^{-1}s\in (UU^{-1})^2$, and similarly $x^{-1}g^{-1}s\in (UU^{-1})^2$. Then $s^{-1}gf^{-1}s=(x^{-1}g^{-1}s)^{-1}(x^{-1}f^{-1}s)\in(UU^{-1})^4$, which contradicts that $gf^{-1}\in FF^{-1}\ssm\{1_\Gamma\}$ and that $(UU^{-1})^4$ is disjoint from the conjugacy classes of $FF^{-1}\ssm\{1_\Gamma\}$. Thus such $f,g\in F$ cannot exist. 

Finally assume $V=xUU^{-1}\in\mathcal{B}$ and $W=yUU^{-1}\in\mathcal{B}$. By a similar argument as the above, we see that if both $FV$ and $FW$ intersect a same member $A\in\mathcal{B}$, then we must have $y\in S\cap F^{-1}Fx(UU^{-1})^4$. For each fixed $g\in F^{-1}F$, the collection of all $y\in S\cap gx(UU^{-1})^4$ gives a collection of left translates $(gx)^{-1}yU$ of $U$ which are pairwise-disjoint subsets of $(UU^{-1})^4U$, and there are at most $M$ such $y$. Since there are at most $\tabs{F}^2$ many such $g\in F^{-1}F$, we see that there are at most $M\cdot\tabs{F}^2$ many $W=yUU^{-1}\in\mathcal{B}$ for which both $FV$ and $FW$ intersect a same member of $\mathcal{B}$. {}
\end{proof}

\begin{lemma}\label{lem:2.11}
Assume $\Gamma$ is a group, $F\subseteq\Gamma$ is a finite set, and $\mathcal{B}$ is a covering of $\Gamma$ satisfying the conclusion of Lemma \ref{lem:2.10}. Let $d,k,n<\omega$ be arbitrary, and assume $F=F_0\cup \ldots\cup F_{n-1}$ where $\tabs{F_0}=\ldots=\tabs{F_{n-1}}=N$ is large enough so that $3kMn^3N^{d+2}(1-k^{-M})^{N-d}\le 1$. Then there is a function $c:\mathcal{B}\to\{0,\ldots,k-1\}$ such that for every $\gamma\in \Gamma$, $i<n$, and $j<k$, there are $>d$ many points in the right translate $F_i\cdot\gamma$ that are only covered by sets $B\in\mathcal{B}$ for which $c(B)=j$. 
\end{lemma}

\begin{proof}
The proof strategy is to invoke the Lov\'{a}sz local lemma. Towards that goal, we will first pick $c:\mathcal{B}\to\{0,\ldots,k-1\}$ independently and uniformly randomly at each point $B\in\mathcal{B}$. We will define a family $\mathcal{A}$ of bad events we wish $c$ to avoid, and we will check that the probability upper-bound $p(\mathcal{A})$ and the dependency degree $d(\mathcal{A})$ satisfy $ep(\mathcal{A})(d(\mathcal{A})+1)\le 1$. Each bad event $A\in\mathcal{A}$ will only depend on finitely many values of $c$, which means it is clopen in the compact product space $k^{\mathcal{B}}$. So a compactness argument implies that the Lov\'{a}sz local lemma holds for $\mathcal{A}$, and there exists a function $c$ avoiding the bad events in $\mathcal{A}$. Finally, we will show that such a function $c$ is as we desired. 

For a set $X\subseteq\Gamma$, define $\mathcal{B}\restrict{}X=\{B\in\mathcal{B}:B\cap X\ne\varnothing\}$. For each $B\in\mathcal{B}$, we define the event $A(B)$ over the domain $\mathcal{B}\restrict{} F\cdot B$ of size $\le M\cdot\tabs{F}$, where $c\in A(B)$ if for some $i<n$ and $j<k$, there are only $d$ or less $f\in F_i$ for which $c[\mathcal{B}\restrict{}fB]=\{j\}$. Let $\mathcal{A}=\{A(B):B\in\mathcal{B}\}$.

We give an upper bound of the probability $p(\mathcal{A})$. Note that if $c\in A(B)$ for some $B\in\mathcal{B}$, then there exists $i<n$, $j<k$, and some size-$d$ subset $D\subseteq F_i$, such that all elements $f\in F_i$ for which $c[\mathcal{B}\restrict{}fB]=\{j\}$ are contained within $D$. For each $f\in F_i\ssm D$, since $\tabs{\mathcal{B}\restrict{}fB}\le M$, the probability that $c$ does not evaluate to $j$ over all of $\mathcal{B}\restrict{}fB$ is at most $1-k^{-M}$. Thus as our assumption on $\mathcal{B}$ implies that $\mathcal{B}\restrict{}fB$ for all $f\in F$ are pairwise-disjoint subsets of $\mathcal{B}$, the probability for a fixed $D$ that all of $f\in F_i\ssm D$ satisfy $c[\mathcal{B}\restrict{}fB]\ne\{j\}$ is at most $(1-k^{-M})^{\tabs{F_i\ssm D}}$. By the union bound, we have the following:
\[
\mathbb{P}[c\in A(B)]\le\sum_{i<n}\sum_{j<k}\binom{\tabs{F_i}}{d}\left(1-k^{-M}\right)^{\tabs{F_i}-d}\le kn\cdot N^d\left(1-k^{-M}\right)^{N-d}
\]
Thus we see that $p(\mathcal{A})\le knN^d(1-k^{-M})^{N-d}$. 

Next, note that two events $A(B)$ and $A(B')$ have intersecting domains exactly when $FB$ and $FB'$ intersect a same member of $\mathcal{B}$. By our assumption on $\mathcal{B}$, every $B$ has at most $M\cdot\tabs{F}^2$ many such $B'$ including $B$ itself, and so we see $d(\mathcal{A})+1\le M\cdot\tabs{F}^2\le Mn^2N^2$. 

By our bounds on $p(\mathcal{A}),d(\mathcal{A})$, we find that we have made $N$ large enough specifically such that $ep(\mathcal{A})(d(\mathcal{A})+1)\le 1$, and the Lov\'{a}sz local lemma applies to give a function $c$ not lying in any bad event in $\mathcal{A}$. 

This $c$ is as we desired, because for every $\gamma\in \Gamma$, if $\gamma\in B\in\mathcal{B}$, then each point $f\gamma\in F\gamma$ is only covered by sets in $\mathcal{B}\restrict{}fB$, and so $c\notin A(B)$ implies that for all $i<n$ and $j<k$, there are $>d$ points $f\gamma$ in $F_i\gamma$ which are covered by sets in $\mathcal{B}\restrict{}fB$ for which $c[\mathcal{B}\restrict{}fB]=\{j\}$. {}
\end{proof}

\begin{theorem}\label{thm:2.12}
Let $\Gamma$ be a locally compact Polish group with a two-sided invariant metric and finite topological dimension. For every $k,n<\omega$, there exists some $N=N(k,n)<\omega$, such that for any sets $F_0,\ldots,F_{n-1}\subseteq\Gamma$ with $\tabs{F_0}=\ldots=\tabs{F_{n-1}}=N$, there exists a sequence of pairwise-disjoint open subsets $D_0,\ldots,D_{k-1}$ of $\Gamma$, for which every right translate $F_i\cdot \gamma$ of every $F_i$ intersects every set $D_j$. In particular, the sequence $\langle D_j:j<k\rangle$ is a domatic $k$-partition with open parts for each of the graphs $\Sch(\Gamma,F_0,\Gamma),\ldots,\Sch(\Gamma,F_{n-1},\Gamma)$. 
\end{theorem}

\begin{proof}
Fix an open cover $\mathcal{B}$ as in Lemma \ref{lem:2.10}, which applies due to Lemmas \ref{lem:2.5} and \ref{lem:2.6}. Fix a function $c:\mathcal{B}\to\{0,\ldots,k-1\}$ as in Lemma \ref{lem:2.11} and a disjoint open family $\mathcal{R}$ as in Theorem \ref{thm:2.9}. Inside every right translate $F_i\cdot \gamma$, apart from the $d=\dim(\Gamma)$ or less points not covered by $\mathcal{R}$ by Theorem \ref{thm:2.9}, Lemma \ref{lem:2.11} guarantees at least one more point covered by $\mathcal{R}$ and colored solely in $j$ by $c$ for any fixed color $j<k$, and in particular its $\mathcal{R}$-cover's associated $\mathcal{B}$-cover (in the sense of Theorem \ref{thm:2.9}) is colored in $j$. Letting $D_j$ be the union of all sets $R\in\mathcal{R}$ whose associated cover in $\mathcal{B}$ is colored in $j$ by $c$ finishes the proof. {}
\end{proof}

\subsection{The open pair property and domatic partitions}\label{sec:2.5}

\begin{definition}\label{def:2.13}
Assuming $\Gamma$ is a group and $D,P\subseteq\Gamma$, we say $D$ \emph{dominates} $P$ if for every $\gamma\in \Gamma$, we have $P\cdot\gamma\cap D\ne\varnothing$. We say an infinite compact Polish group $\Gamma$ has the \emph{open pair property} if for every finite collection of nonempty perfect subsets $P_0,\ldots,P_{n-1}$ of $\Gamma$, there exists a pair of disjoint open sets $A_0,A_1\subseteq\Gamma$ each of which dominates all of $P_0,\ldots,P_{n-1}$. 
\end{definition}

\begin{lemma}\label{lem:2.14}
Let $\Gamma$ be an infinite compact Polish group with finite topological dimension. Then $\Gamma$ has the open pair property. 
\end{lemma}

\begin{proof}
Since every compact Polish group admits a two-sided invariant metric \citep[Exercise 2.1.5]{gao}, Theorem \ref{thm:2.12} applies. Let $P_0,\ldots,P_{n-1}$ be nonempty perfect subsets of $\Gamma$, and for each $i<n$ let $F_i\subseteq P_i$ be any finite subset of size $N=N(2,n)$ as in Theorem \ref{thm:2.12}. Then the theorem gives a pair of disjoint open sets $D_0,D_1$ each of which dominates every $F_i$, and we see that they also dominate every $P_i$ because $F_i\subseteq P_i$. {}
\end{proof}

\begin{lemma}\label{lem:2.15}
Let $\Gamma$ be a compact Polish group, and assume $D\subseteq \Gamma$ is an open set that dominates some set $P\subseteq\Gamma$. Then there exists a finite subset $F\subseteq P$ such that $D$ dominates $F$. Moreover, if $d$ is a compatible two-sided invariant metric on $\Gamma$, then there exists some $r>0$ such that for every $\gamma\in \Gamma$ there exists $f\in F$ for which $B_d(f,r)\cdot \gamma\subseteq D$. It follows that there is an open set $U$ such that $U\subseteq\overline{U}\subseteq D$ and $U$ dominates both $F$ and $P$. 
\end{lemma}

\begin{proof}
First, note that $D$ dominates $P$ if and only if for every $\gamma\in \Gamma$ there exists $p\in P$ such that $p\cdot\gamma\in D$, which is also equivalent to $\{p^{-1}D:p\in P\}$ being an open cover of $\Gamma$. So since $\Gamma$ is compact, there is a finite subcover $\{p^{-1}D:p\in F\}$ for some finite $F\subseteq P$, and we see that $D$ dominates this finite subset $F$ of $P$. 

For the next part of the lemma, fix a compatible two-sided invariant metric $d$ on $\Gamma$. Define the function $g:\Gamma\to \R\cup\{\infty\}$ via $g(\gamma)=\max_{f\in F}d(f\cdot \gamma,\Gamma\ssm D)$. Note that for each $f\in F$, the function $\gamma\mapsto d(f\cdot \gamma,\Gamma\ssm D)$ is continuous in $\gamma$, and so $g$ as a finite maximum of continuous functions is also continuous. For each $\gamma$, by assumption there exists some $f\in F$ such that $f\cdot \gamma\in D$, which means $g(\gamma)\ge d(f\cdot\gamma,\Gamma\ssm D)>0$ for this $f$. Consequently $g$ is a continuous function over a compact domain $\Gamma$ whose range is within $(0,\infty]$, and so compactness allows us to put a positive lower-bound $r>0$ on $g[\Gamma]$. 

By definition of $g$, we see that for all $\gamma\in \Gamma$ there exists some $f\in F$ such that $B_d(f\cdot \gamma,r)\subseteq D$. Right-invariance of $d$ gives us $B_d(f\cdot\gamma,r)=B_d(f,r)\cdot\gamma\subseteq D$. 

Finally, let $U=\{x\in \Gamma:d(x,\Gamma\ssm D)>r/2\}$, so that $U\subseteq\overline{U}\subseteq D$. For every $\gamma\in \Gamma$ there is some $f\in F$ such that $B_d(f\cdot \gamma,r)\subseteq D$, which means that $f\cdot\gamma\in U$, which means that $U$ dominates $F$. Since $F\subseteq P$, we see $U$ also dominates $P$. {}
\end{proof}

\begin{lemma}\label{lem:2.16}
Let $\Gamma$ be a compact Polish group with the open pair property. Assume that $U\subseteq\Gamma$ is an open set which dominates finitely many nonempty perfect sets $P_0,\ldots,P_{n-1}\subseteq\Gamma$. Then there exists a pair of disjoint open subsets $A_0,A_1\subseteq U$ each of which dominates all of $P_0,\ldots,P_{n-1}$. 
\end{lemma}

\begin{proof}
Fix a compatible two-sided invariant metric $d$ on $\Gamma$. By Lemma \ref{lem:2.15}, for each $i<n$, there is a finite set $F_i\subseteq P_i$ and $r>0$, such that every $\gamma\in \Gamma$ has some $f\in F_i$ for which $B_d(f,r)\cdot\gamma\subseteq U$. For each $f\in F_i\subseteq P_i$, we fix a nonempty perfect subset $P_{i,f}\subseteq P_i\cap B_d(f,r)$, such that every $\gamma\in \Gamma$ has some $f\in F_i$ for which $P_{i,f}\cdot\gamma\subseteq B_d(f,r)\cdot\gamma\subseteq U$. 

Let $\mathcal{P}=\{P_{i,f}:i<n,f\in F_i\}$, which is a finite collection of nonempty perfect sets. By the open pair property of $\Gamma$, there are disjoint open sets $D_0,D_1\subseteq \Gamma$ each of which dominates all of $\mathcal{P}$. Let $A_0=D_0\cap U$ and $A_1=D_1\cap U$, which we will show are as desired. 

Fix $i<n$, $j\in\{0,1\}$, and $\gamma\in \Gamma$, and it remains to show that $P_i\cdot\gamma\cap A_j\ne\varnothing$. By the above, there exists some $f\in F_i$ for which $P_{i,f}\cdot\gamma\subseteq U$. Since $D_j$ dominates $P_{i,f}\in\mathcal{P}$, there is some $p\in P_{i,f}$ such that $p\cdot\gamma\in D_j$. Then $p\cdot\gamma\in D_j\cap U=A_j$ and also $p\in P_{i,f}\subseteq P_i$, which means that $p\cdot\gamma\in P_i\cdot\gamma\cap A_j\ne\varnothing$ as desired. {}
\end{proof}

\begin{theorem}\label{thm:2.17}
Let $\Gamma$ be a finite-dimensional compact Polish group. Let $S_0,\ldots,S_{n-1}\subseteq\Gamma$ be subsets such that every $\overline{S_i}\subseteq\Gamma$ is uncountable. Then there is a sequence of pairwise-disjoint open sets $\langle D_j:j<\omega\rangle$, for which every right translate $S_i\cdot\gamma$ of every $S_i$ intersects every set $D_j$. In particular, the sequence $\langle D_j:j<\omega\rangle$ is a domatic $\aleph_0$-partition with open parts for each of the graphs $\Sch(\Gamma,S_0,\Gamma),\ldots,\Sch(\Gamma,S_{n-1},\Gamma)$. 
\end{theorem}

\begin{proof}
By the perfect set theorem, for each $i<n$ fix a nonempty perfect set $P_i\subseteq\overline{S_i}$. 

We first construct $\langle D_j:j<\omega\rangle$ which dominates every $P_i$. Note that $U_0=\Gamma$ is open and dominates every $P_i$. Inductively, whenever $\{D_0,\ldots,D_{n-1},U_n\}$ is a pairwise-disjoint open family which dominates every $P_i$, we can apply Lemmas \ref{lem:2.14} and \ref{lem:2.16} to split $U_n$ into disjoint open parts $D_n,U_{n+1}\subseteq U_n$ which still dominate every $P_i$, and thus we get a pairwise-disjoint open family $\{D_0,\ldots,D_n,U_{n+1}\}$ which dominates every $P_i$. Carrying on this induction through $n<\omega$, we get a pairwise-disjoint open family $\{D_j:j<\omega\}$ that dominates every $P_i$. 

It remains to show every $D_j$ also dominates every $S_i$. For every $\gamma\in \Gamma$, we have $P_i\gamma\cap D_j\ne\varnothing$, and every point in $P_i \gamma$ is a limit point of $S_i\gamma$. So since $D_j$ is open and contains some limit point of $S_i\gamma$, it also contains some point of $S_i\gamma$. We get that $S_i\gamma\cap D_j\ne\varnothing$, which means $D_j$ also dominates $S_i$ as desired. {}
\end{proof}

\begin{corollary}\label{cor:2.18}
Let $\Gamma$ be a finite-dimensional compact Polish group, and let $S\subseteq\Gamma$ be a subset. Then the graph $\Sch(\Gamma,S,\Gamma)$ admits a domatic $\aleph_0$-partition with open parts, if and only if it admits a domatic $\aleph_0$-partition with Baire measurable parts, if and only if $\overline{S}\subseteq\Gamma$ is uncountable. 
\end{corollary}

\begin{proof}
Follows directly from Theorems \ref{thm:2.1} and \ref{thm:2.17}. {}
\end{proof}

\begin{corollary}\label{cor:2.19}
Let $n\in\N$, and let $S\subseteq\R^n$ be a subset of $\R^n$. Then the graph $\Sch(\R^n,S,\R^n)$ admits a domatic $\aleph_0$-partition with open or Baire measurable parts if and only if either $\overline{S}\subseteq\R^n$ is uncountable or $S\subseteq\R^n$ is unbounded. 
\end{corollary}

\begin{proof}
We split into three cases: The case when $\overline{S}$ is countable and bounded, the case when $\overline{S}$ is uncountable, and the case when $S$ is unbounded. 

When $\overline{S}$ is countable and bounded, Theorem \ref{thm:2.1} implies that $\Sch(\R^n,S,\R^n)$ does not admit a domatic $\aleph_0$-partition with Baire measurable parts. 

When $\overline{S}$ is uncountable, the set $\overline{\pi(S)}\subseteq\R^n/\Z^n$ is also uncountable, where $\pi:\R^n\to \R^n/\Z^n$ is the usual projection map, and this is because $\overline{\pi(S)}\supseteq\pi(\overline{S})$ and $\pi$ is countable-to-one. Then Theorem \ref{thm:2.17} gives a continuous domatic partial function $f:\R^n/\Z^n\partialto{}\aleph_0$ for the graph $\Sch(\R^n/\Z^n,\pi(S),\R^n/\Z^n)$ on $\R^n/\Z^n$, and one can check that the pullback $f\circ \pi:\R^n\partialto{}\aleph_0$ is a continuous domatic partial function for $\Sch(\R^n,S,\R^n)$. 

When $S$ is unbounded, we may use diagonalization to get a rapidly increasing sequence of positive radii $0= R_0<R_1<\ldots$ such that for any $i<\omega$ and any $x\in\R^n$ with $\|x\|<i+1$, the translate $S+x$ intersects the open spherical shell $B(0,R_{i+1})\ssm\overline{B(0,R_i)}$. Thus any translate of $S$ eventually intersects with every far enough open spherical shell in this sequence. Then an $\aleph_0$-coloring of these spherical shells where each color is used infinitely often will give a continuous domatic partial function $f:\R^n\partialto{}\aleph_0$ for $\Sch(\R^n,S,\R^n)$. {}
\end{proof}

Curiously, even among the small class of compact Polish groups, our analysis above leaves open the case of infinite dimension. For now, the following question remains open:

\begin{question}\label{que:2.20}
Let $P$ be a nonempty perfect subset of the Polish group $\Gamma=(\R/\Z)^\omega$. Does the graph $G=\Sch(\Gamma,P,\Gamma)$ admit domatic $\aleph_0$-partitions with open, Borel, or Baire measurable parts? Does $G$ even admit a domatic bipartition with open parts? 
\end{question}

See also Corollary \ref{cor:3.6} for the existence of measure-theoretic domatic $\aleph_0$-partitions on the Schreier graph $\Sch(\Gamma,S,\Gamma)$. 

\subsection{Application to sum sets}\label{sec:2.6}

In this section, we use ideas from the previous sections to prove Theorem \ref{thm:2.27} and its Corollary \ref{cor:2.29}, which extend a result on sum sets by Erd\H{o}s--Kunen--Mauldin \citep[Theorem 1]{ekm}. 

\begin{definition}
Assume $D,P\subseteq\Gamma$ are subsets of a group $\Gamma$, we say that $D$ \emph{additively dominates} $P$ if $P\cdot D=\Gamma$. When $D^2=\Gamma$, we say that $D$ \emph{additively dominates itself}. 
\end{definition}

Recall from the proof of Lemma \ref{lem:2.15} that $D$ dominates $P$ if and only if $P^{-1}\cdot D=\Gamma$. This means that $D$ additively dominates $P$ if and only if $D$ dominates $P^{-1}$. 

\begin{theorem}\label{thm:2.22}
Let $\Gamma$ be a finite-dimensional compact Polish group, and let $P\subseteq \Gamma$ be a nonempty closed perfect subset. Then there exists a family $\langle D_i:i<2^{\aleph_0}\rangle$ of $2^{\aleph_0}$ pairwise-disjoint closed subsets of $\Gamma$, each of which additively dominates $P$. 
\end{theorem}

\begin{proof}
By Lemma \ref{lem:2.16}, we can build a tree $\{ U_s:s\in {}^{\omega>}2\}$ of open subsets of $\Gamma$, such that $U_{\langle\:\rangle}=\Gamma$, for each $s\preceq t$ we have $U_s\supseteq U_t$, for each $s$ we have $U_{s\concat{}0}\cap U_{s\concat{}1}=\varnothing$, and every $U_s$ dominates $P^{-1}$. By the last part of Lemma \ref{lem:2.15}, we may also shrink each $U_s$ along the construction in such a way that $\overline{U_{s\concat{}0}},\overline{U_{s\concat{}1}}\subseteq U_s$ for every $s$. 

For every $x\in {}^\omega 2$, let $D_x=\bigcap_{n<\omega}U_{x\restrict{} n}=\bigcap_{n<\omega}\overline{U_{x\restrict{} n}}$, and we claim that $D_x$ dominates $P^{-1}$. Fix a $\gamma\in \Gamma$, and we'll show $P^{-1}\gamma\cap D_x\ne\varnothing$. For every $n<\omega$, since $U_{x\restrict{}n}$ dominates $P^{-1}$, there exists some $p_n\in P^{-1}$ such that $p_n\gamma\in U_{x\restrict{}n}$. Since $P^{-1}\subseteq\Gamma$ is compact, there is a subsequence of $p_n$'s which converges to some $p\in P^{-1}$. Since $\{U_{x\restrict{}n}:n<\omega\}$ is a decreasing family, taking the limit of $p_n\gamma\in U_{x\restrict{}n}$ along this subsequence gives $p\gamma\in\bigcap_{n<\omega}\overline{U_{x\restrict{}n}}=D_x$, which means $D_x$ dominates $P^{-1}$. Thus $\{D_x:x\in {}^\omega 2\}$ is the family that we wanted. {}
\end{proof}

\begin{lemma}\label{lem:2.23}
Let $\Gamma$ be a non-Boolean connected Polish group, where a group is Boolean if every non-identity element has order $2$. For every $x\in \Gamma$, the closed set $\{\gamma\in \Gamma:\gamma^2=x\}$ is nowhere dense. In particular, every nontrivial connected locally compact Polish group has this property. 
\end{lemma}

\begin{proof}
We first assume that $\Gamma$ is a connected Polish group and $x\in\Gamma$ is such that $\{\gamma\in \Gamma:\gamma^2=x\}$ has nonempty interior. We will show that such a $\Gamma$ is Boolean. 

Fix a nonempty open set $U\subseteq\Gamma$ such that $\gamma^2=x$ for all $\gamma\in U$. Let $g\in U$, and let $V\ni 1_\Gamma$ be a symmetric open neighborhood of the identity such that $V^2g\subseteq U$. For every $h\in V^2$, we have $hg,g\in U$, which means that $ghg^{-1}=h^{-1}(hghg)(g^{-1}g^{-1})=h^{-1}xx^{-1}=h^{-1}$. 

Assume $a,b\in V$, so then $a,b,ab\in V^2$, and the above argument implies that $gag^{-1}=a^{-1}$, $gbg^{-1}=b^{-1}$, and $a^{-1}b^{-1}=(gag^{-1})(gbg^{-1})=g(ab)g^{-1}=(ab)^{-1}=b^{-1}a^{-1}$. We find that $a^{-1}$ commutes with $b^{-1}$ for all $a,b\in V$, and since $V$ is symmetric, $a$ and $b$ also commute for every $a,b\in V$. 

Since $\Gamma$ is connected and $V$ is an open neighborhood of $1_\Gamma$, $V$ is a generating set of $\Gamma$. Since $\Gamma$ has a commutative generating set, $\Gamma$ is abelian. This then means for all $h\in V\subseteq V^2$, we have $h=ghg^{-1}=h^{-1}$, and so every $h\in V$ must have order at most $2$. Finally since $\Gamma$ is an abelian group generated by a set $V$ of elements of order at most $2$, $\Gamma$ is Boolean as well. 

For the next part of the lemma, we will show that every nontrivial connected locally compact Polish group is non-Boolean. By Corollary \ref{cor:2.3}, every such group is an inverse limit of nontrivial connected Lie groups. We may notice that every nontrivial connected Lie group is non-Boolean, by considering group elements near the identity. Thus their inverse limit must also be non-Boolean, which finishes the proof. {}
\end{proof}

\begin{lemma}\label{lem:2.24}
Let $\Gamma$ be a perfect Polish group, and assume $D\subseteq\Gamma$ is an open set that dominates some set $P\subseteq\Gamma$, where $P$ has no isolated points. Then for any finite subset $F\subseteq D$, the set $D\ssm F$ also dominates $P$. 
\end{lemma}

\begin{proof}
For every $\gamma\in \Gamma$, the set $P\gamma\cap D\ne\varnothing$ is nonempty and has no isolated points. So $P\gamma\cap D$ is infinite, and $F$ being finite implies $P\gamma\cap (D\ssm F)=(P\gamma\cap D)\ssm F\ne\varnothing$ is nonempty, which means $D\ssm F$ dominates $P$. {}
\end{proof}

\begin{lemma}\label{lem:2.25}
Let $\Gamma$ be an infinite connected compact Polish group, and assume $U\subseteq\Gamma$ is an open set that additively dominates itself. Then there exists a finite subset $F\subseteq U$ such that $U\ssm F$ additively dominates $F$. 
\end{lemma}

\begin{proof}
Since $U$ dominates $U^{-1}$, Lemma \ref{lem:2.15} gives a finite subset $F_0^{-1}\subseteq U^{-1}$ such that $U$ dominates $F_0^{-1}$. By Lemma \ref{lem:2.24}, $U\ssm F_0$ dominates $U^{-1}$, and by another application of Lemma \ref{lem:2.15}, there is a finite collection $B_0^{-1},\ldots,B_{n-1}^{-1}$ of nonempty open subsets of $U^{-1}$, such that for every $\gamma\in \Gamma$, there exists some $i<n$ for which $B_i^{-1}\gamma\subseteq U\ssm F_0$. It follows that for any finite tuple $(f_0^{-1},\ldots,f_{n-1}^{-1})\in B_0^{-1}\times\ldots\times B_{n-1}^{-1}$, the set $U\ssm F_0$ dominates $\{f_0^{-1},\ldots,f_{n-1}^{-1}\}\subseteq U^{-1}$. 

We claim that for a comeager set of tuples $(f_0,\ldots,f_{n-1})\in \Gamma^n$, the finite set $F_1=\{f_0,\ldots,f_{n-1}\}$ satisfies $F_1F_1F^{-1}_1\cap F_0=\varnothing$ and $F_1F_1\cap F_0F_0=\varnothing$. 

To fulfill the first condition $F_1F_1F_1^{-1}\cap F_0=\varnothing$, it suffices to guarantee for every $(i,j,k)\in\{0,\ldots,n-1\}^3$ and $g\in F_0$, that $f_if_jf_k^{-1}\ne g$ happens comeagerly often in $\Gamma^n$. If one of $i,j,k$ is different from the other two, then $f_if_jf_k^{-1}\ne g$ happens comeagerly if we fix the other two $f$'s and move the one different $f$ freely, and so $f_if_jf_k^{-1}\ne g$ happens comeagerly over $\Gamma^n$ by Kuratowski--Ulam \citep[Theorem 8.41]{kechris}. Otherwise if $i=j=k$, then $f_i=f_if_jf_k^{-1}\ne g$ also happens comeagerly over $\Gamma^n$. 

To fulfill the second condition $F_1F_1\cap F_0F_0=\varnothing$, one can apply a similar argument as the previous case, using the fact that Lemma \ref{lem:2.23} applies. 

Therefore since $B_0\times\ldots\times B_{n-1}\subseteq\Gamma^n$ is nonmeager, we can fix some tuple $(f_0,\ldots,f_{n-1})\in B_0\times\ldots\times B_{n-1}$ such that the finite set $F_1=\{f_0,\ldots,f_{n-1}\}$ satisfies $F_1F_1F_1^{-1}\cap F_0=\varnothing$ and $F_1F_1\cap F_0F_0=\varnothing$. We note here that $F_0F_1\cap F_1F_1=\varnothing$, since otherwise there exists $(g,f)\in F_0\times F_1$ where $gf\in F_1F_1$, and then $g\in F_1F_1F_1^{-1}\cap F_0\ne\varnothing$ gives a contradiction. In other words, we now know that $(F_0F_0\cup F_0F_1)\cap F_1F_1=\varnothing$. 

Let $F=F_0\cup F_1\subseteq U$, and recall from earlier that $U$ dominates $F_0^{-1}$ and $U\ssm F_0$ dominates $F_1^{-1}$. This means that $F_0U=F_1(U\ssm F_0)=\Gamma$. We see that the following holds:
\begin{align*}
F\cdot(U\ssm F)&=(F_0\cup F_1)\cdot\left(U\ssm (F_0\cup F_1)\right)\\
&=F_0\cdot\left(U\ssm (F_0\cup F_1)\right)\cup F_1\cdot\left((U\ssm F_0)\ssm F_1\right)\\
&\supseteq \left(\Gamma\ssm (F_0F_0\cup F_0F_1)\right)\cup\left(\Gamma\ssm F_1F_1\right)\\
&=\Gamma\ssm \left((F_0F_0\cup F_0F_1)\cap F_1F_1\right)\\
&=\Gamma
\end{align*}
We conclude that $U\ssm F$ additively dominates $F$. {}
\end{proof}

\begin{lemma}\label{lem:2.26}
Let $\Gamma$ be an infinite connected compact Polish group with the open pair property. Assume that $U\subseteq\Gamma$ is an open set which additively dominates finitely many perfect sets $P_0,\ldots,P_{m-1}$ and $U$ itself, and assume that finitely many open sets $Q_0,\ldots,Q_{n-1}$ additively dominate $U$. Then there exists a pair of disjoint open subsets $A_0,A_1\subseteq U$, each of which additively dominates all of $A_0,A_1,P_0,\ldots,P_{m-1}$ and is additively dominated by $Q_0,\ldots,Q_{n-1}$, and moreover $\overline{A_0},\overline{A_1}\subseteq U$. 
\end{lemma}

\begin{proof}
By Lemma \ref{lem:2.15} and since domination is closed upwards, there is a finite subset $F_Q\subseteq U$ such that every $Q_j$ additively dominates $F_Q$. By Lemma \ref{lem:2.24}, $U\ssm F_Q$ additively dominates $U$, and a similar argument shows that $U\ssm F_{Q}$ additively dominates itself. By Lemma \ref{lem:2.25}, there is a finite subset $F_U\subseteq U\ssm F_Q$ such that $U\ssm (F_Q\cup F_U)$ additively dominates $F_U$. Letting $F=F_Q\cup F_U\subseteq U$, we see that by upward closure of domination, the sets $Q_0,\ldots,Q_{n-1},U\ssm F$ all additively dominate $F$, and by Lemma \ref{lem:2.24}, the set $U\ssm F$ additively dominates all of the sets $P_0,\ldots,P_{m-1}$. 

By Lemma \ref{lem:2.15}, there is an open subset $V\subseteq\overline{V}\subseteq U\ssm F$ such that $V$ additively dominates all of the sets $F,P_0,\ldots,P_{m-1}$, which is because we can take a finite union of all such $V$'s over each of the sets $F,P_0,\ldots,P_{m-1}$.  Also by Lemma \ref{lem:2.15}, if $d$ is a compatible two-sided invariant metric on $\Gamma$, then for every sufficiently small $r>0$, every set that intersects every ball $B_d(f,r)$ for $f\in F$ is additively dominated by the sets $Q_0,\ldots,Q_{n-1},V$ which additively dominated $F$, and the radius-$r$ ball $B_d(F,r)$ around $F$ satisfies $\overline{B_d(F,r)}\subseteq U\ssm \overline{V}$. 

Since $\Gamma$ is perfect, we can take a pair of disjoint open subsets $W_0,W_1\subseteq B_d(F,r)$, such that each of $W_0,W_1$ intersects every ball $B_d(f,r)$ for $f\in F$. Thus each of $W_0,W_1$ is additively dominated by the sets $Q_0,\ldots,Q_{n-1},V$. 

The open set $V$ additively dominates the perfect sets $P_0,\ldots,P_{m-1},\overline{W_0},\overline{W_1}$, and so by Lemma \ref{lem:2.16}, there are disjoint open subsets $D_0,D_1\subseteq V$ each of which additively dominates the sets $P_0,\ldots,P_{m-1},\overline{W_0},\overline{W_1}$. By an argument in the proof of Theorem \ref{thm:2.17}, $D_0,D_1$ each also additively dominates the sets $W_0,W_1$. 

Let $A_0=D_0\cup W_0$ and $A_1=D_1\cup W_1$. Then $A_0,A_1$ additively dominate each other and themselves since each of $D_0,D_1$ additively dominates each of $W_0,W_1$, also $A_0,A_1$ additively dominate $P_0,\ldots,P_{m-1}$ since $D_0,D_1$ do, $A_0,A_1$ are additively dominated by $Q_0,\ldots,Q_{n-1}$ since $W_0,W_1$ are, and $\overline{A_0},\overline{A_1}\subseteq\overline{B_d(F,r)}\cup \overline{V}\subseteq U$. {}
\end{proof}

\begin{theorem}\label{thm:2.27}
Let $\Gamma$ be an infinite finite-dimensional connected compact Polish group, and let $P\subseteq\Gamma$ be a nonempty closed perfect subset. Then there exists a family $\langle D_i:i<2^{\aleph_0}\rangle$ of $2^{\aleph_0}$ pairwise-disjoint closed subsets of $\Gamma$, such that every $D_i$ additively dominates $P$, and every $D_i$ additively dominates every $D_j$, for all $i,j<2^{\aleph_0}$. 
\end{theorem}

\begin{proof}
Like Theorem \ref{thm:2.22}, we can build a tree $\{U_s:s\in {}^{\omega>}2\}$ of open subsets of $\Gamma$, such that $U_{\langle\:\rangle}=\Gamma$, for each $s\prec t$ we have $U_s\supseteq \overline{U_t}$, for each $s$ we have $U_{s\concat{}0}\cap U_{s\concat{}1}=\varnothing$, every $U_s$ additively dominates $P$, and every $U_s$ additively dominates every $U_t$, for all $s,t\in {}^{\omega>}2$. 

The strategy is to build $\{U_s:s\in {}^{\omega>}2\}$ ``one-by-one'', where we visit every node in the tree ${}^{\omega>}2$ in $\omega$ stages, while increasingly traversing along each branch. At a given stage $s\in{}^{\omega>}2$, we split the open set $U_s$ into two disjoint open subsets $U_{s\concat{}0},U_{s\concat{}1}$ using Lemma \ref{lem:2.26}, such that among the open sets $U_t$ we've constructed so far, each one additively dominates $P$ and all of them additively dominate each other and themselves. In using Lemma \ref{lem:2.26}, note that $U_s$ additively dominates $U_t$ if and only if $U_s$ additively dominates the perfect set $\overline{U_t}$ by an argument in the proof of Theorem \ref{thm:2.17}. 

For every $x\in {}^\omega 2$, let $D_x=\bigcap_{n<\omega}U_{x\restrict{} n}=\bigcap_{n<\omega}\overline{U_{x\restrict{} n}}$. A same argument as Theorem \ref{thm:2.22} implies that the family $\{D_x:x\in {}^\omega 2\}$ is as we desired. {}
\end{proof}

\begin{theorem}\label{thm:2.28}
Let $\Gamma$ be an infinite connected compact Polish group. Then there exists a family $\langle D_i:i<2^{\aleph_0}\rangle$ of $2^{\aleph_0}$ pairwise-disjoint closed subsets of $\Gamma$, such that every $D_i$ additively dominates every $D_j$, for all $i,j<2^{\aleph_0}$. 
\end{theorem}

\begin{proof}
We first prove a weak version of Lemma \ref{lem:2.14} without the finite dimension assumption. Namely, let $\Gamma$ be a perfect Polish group, and let $P_0,\ldots,P_{n-1}\subseteq \Gamma$ be nonempty open sets. Then we will show that there exists a pair of disjoint open sets $A_0,A_1\subseteq\Gamma$, each of which dominates the sets $P_0,\ldots,P_{n-1}$. 

By Birkhoff--Kakutani \citep[Theorem 2.1.1]{gao}, $\Gamma$ has a compatible right-invariant metric $d$. Then there is some $r>0$ such that every right translate $P_i\gamma$ of a set $P_i$ contains some ball of $d$-radius-$r$. By a ball-packing argument and since $\Gamma$ is perfect, there are disjoint discrete sets $S_0,S_1\subseteq \Gamma$ which both intersect every radius-$r$ ball in $\Gamma$. Then let $A_0,A_1$ be disjoint open sets which separate $S_0,S_1$, and our claimed result follows. 

Then by repeating the proofs of Lemma \ref{lem:2.26} and Theorem \ref{thm:2.27} with the above notion of a weaker open pair property in place of the normal one, the full result follows. {}
\end{proof}

\begin{corollary}\label{cor:2.29}
Let $1\le n\in\N$, and let $P\subseteq\R^n$ be a nonempty closed perfect subset of $\R^n$. Then there exists a family $\langle C_i:i<2^{\aleph_0}\rangle$ of $2^{\aleph_0}$ pairwise-disjoint closed subsets of $\R^n$, such that $P+C_i=\R^n$ and $C_i+C_j=\R^n$ for all $i,j<2^{\aleph_0}$. 
\end{corollary}

\begin{proof}
Let $\pi:\R^n\to\R^n/\Z^n$ be the usual projection map, and fix a nonempty closed perfect subset $K\subseteq \pi[P]$ in $\R^n/\Z^n$. By Theorem \ref{thm:2.27}, there is a family $\langle D_i:i<2^{\aleph_0}\rangle$ of pairwise-disjoint closed subsets of $\R^n/\Z^n$ which additively dominate $K$, each other, and themselves. Similar to Corollary \ref{cor:2.19}, the pullbacks $C_i=\pi^{-1}[D_i]$ for $i<2^{\aleph_0}$ are as desired. {}
\end{proof}

\section{More examples of domatic \texorpdfstring{$\aleph_0$}{infinite}-partitions}\label{sec:3}

\subsection{A greedy algorithm on smooth Borel graphs}\label{sec:3.1}

\begin{theorem}\label{thm:3.1}
Let $G$ be an out-degree $\aleph_0$-regular Borel graph with countable in-degrees on a Borel space $X$ of vertices, such that its connectedness (countable Borel) equivalence relation $E_G$ on $X$ is smooth. Then $G$ admits a Borel domatic $\aleph_0$-partition. 
\end{theorem}

\begin{proof}
First we describe a greedy algorithm performed on a connected out-degree $\aleph_0$-regular graph $G$ with a fixed enumeration $V=\{v_0,v_1,v_2,\ldots\}$ of its countable vertex set, which outputs a domatic $\aleph_0$-partition for the countable graph $G$. 

We start with the empty coloring $f_0=\varnothing$ on $V$, and we fix a countable enumeration $\{(v_i,c_i):i<\omega\}=\omega\times \omega$. At stage $i<\omega$ of this algorithm, we extend a finite partial coloring $f_i$ on $V$ to a finite partial coloring $f_{i+1}\supseteq f_i$ on $V$, such that $v_i\in\dom(f_{i+1})$, and $f_{i+1}$ colors the first $f_i$-uncolored neighbor of $v_i$ in color $c_i$, if $f_i$ hasn't colored any neighbor of $v_i$ in color $c_i$ yet. We can check that this algorithm always outputs a full-domain domatic $\aleph_0$-partition for a connected out-degree $\aleph_0$-regular countable graph $G$. 

Using standard descriptive graph combinatorial arguments and Lusin--Novikov uniformization \citep[Theorem 18.10]{kechris}, one sees that the above algorithm can be performed on locally countable Borel graphs when $E_G$ is smooth \citep[Exercise 18.20]{kechris}, so that it outputs Borel colorings. {} 
\end{proof}

\subsection{Measure-theoretic domatic \texorpdfstring{$\aleph_0$}{infinite}-partitions}\label{sec:3.2}

In this section, we prove Theorem \ref{thm:3.5}. 

\begin{lemma}\label{lem:3.2}
Let $(X,\mu)$ be a Borel probability space. Let $G$ be an out-degree $\aleph_0$-regular Borel graph with countable in-degrees on the vertex set $X$, and assume there exists a Borel function $f:X\to\omega$ such that every vertex $x\in X$ has an infinitely colored out-neighborhood, meaning $\tabs{f[N_G(x)]}=\aleph_0$. Then $G$ admits a $\mu$-measurable domatic $\aleph_0$-partition. 
\end{lemma}

\begin{proof}
Let $\kappa$ be the probability measure on $\omega$ defined by $\kappa(\{n\})=2^{-n-1}$, and let $\lambda=\prod_{i<\omega}\kappa$ be the product Borel probability measure on $\prod_{i<\omega}\omega={}^\omega \omega$. Note that since $\kappa$ does not vanish on singletons, it follows that for every infinite subset $A\subseteq \omega$, there exists a $\lambda$-conull set of functions $r:\omega\to\omega$ such that $r[A]=\omega$. 

Let $B=\{(x,r)\in X\times{}^\omega\omega:\forall y\in[x]_{E_G}(r\circ f[N_G(y)]=\omega)\}$. Then since $E_G$ is a countable Borel equivalence relation, $B$ is Borel. Since for every $y\in X$, the set $f[N_G(y)]\subseteq\omega$ is infinite, the previous arguments imply that every section $B_x\subseteq{}^\omega\omega$ of $B$ is $\lambda$-conull. By Fubini's theorem, there exists an $r\in {}^\omega\omega$ such that the section $B^r\subseteq X$ is $\mu$-conull. This means that the Borel function $r\circ f:X\to\omega$ is domatic at every $x\in B^r= [B^r]_{E_G}$. 

Since $B^r$ is $E_G$-invariant, we can modify the function $r\circ f:X\to\omega$ over the $E_G$-invariant $\mu$-null set $X\ssm B^r$ in the same way as Theorem \ref{thm:3.1}, so that we get a $\mu$-measurable function $g:X\to\omega$ domatic everywhere for $G$ as desired. {}
\end{proof}

\begin{lemma}\label{lem:3.3}
Let $(X,\tau)$ be a Polish space. Let $G$ be an out-degree $\aleph_0$-regular Borel graph with countable in-degrees on the vertex set $X$, and assume there exists a Borel function $f:X\to\omega$ such that every vertex $x\in X$ has an infinitely colored out-neighborhood, meaning $\tabs{f[N_G(x)]}=\aleph_0$. Then $G$ admits a $\tau$-Baire measurable domatic $\aleph_0$-partition. 
\end{lemma}

\begin{proof}
The proof is the same as that of Lemma \ref{lem:3.2}, with all measure-theoretic elements replaced by their Baire category counterparts. {}
\end{proof}

\begin{lemma}\label{lem:3.4}
Let $(X,\mu)$ be a Borel probability space. Let $G$ be an out-degree $\aleph_0$-regular Borel graph on the vertex set $X$. Then for every $k<\omega$ and $\varepsilon>0$, there exists a Borel function $f:X\to \{0,\ldots,k-1\}$ that is domatic at a set of vertices of $\mu$-measure $\ge 1-\varepsilon$. 
\end{lemma}

\begin{proof}
Let $n<\omega$ be sufficiently large so that $k(1-k^{-1})^n\le\varepsilon/2$. 

By Lusin--Novikov uniformization \citep[Theorem 18.10]{kechris}, there are Borel functions $f_0,\ldots,f_{n-1}:X\to X$, such that for every $x\in X$, the elements $f_0(x),\ldots,f_{n-1}(x)\in N_G(x)$ are $n$ distinct $G$-neighbors of $x$. 

Let $\{A_0,A_1,\ldots\}$ be a countable family of Borel subsets of $X$ that separates points. For $i<\omega$, let $\mathcal{P}_i$ be the finite Borel partition of $X$ generated by $\{A_0,\ldots,A_{i-1}\}$. Then for every $x\in X$, the $n$ distinct points $f_0(x),\ldots,f_{n-1}(x)$ will be completely separated by some partition $\mathcal{P}_i$ (and every partition that comes afterwards), and so we can pick a sufficiently fine partition $\mathcal{P}=\mathcal{P}_i$ such that the Borel set $B\subseteq X$ of all elements $x\in X$ for which $f_0(x),\ldots,f_{n-1}(x)$ are completely separated by $\mathcal{P}$ satisfies $\mu(B)\ge 1-\varepsilon/2$. 

Since $\mathcal{P}$ is a finite Borel partition of $X$, we can uniformly randomly pick a function $c:\mathcal{P}\to \{0,\ldots,k-1\}$, and let $f_c:X\to \{0,\ldots,k-1\}$ be its corresponding Borel coloring on $X$ that assigns the color $c(P)$ to every element inside $P\in\mathcal{P}$. For every $x\in B$, the probability that $f_c$ is domatic at $x$ is lower-bounded by the probability that $f_c$ colors the $\mathcal{P}$-separated points $f_0(x),\ldots,f_{n-1}(x)$ with all colors in $\{0,\ldots,k-1\}$, which is at least $1-k(1-k^{-1})^n\ge 1-\varepsilon/2$. Thus the product probability measure of the set of all pairs $(c,x)$ for which $f_c$ is domatic at $x$ is at least $(1-\varepsilon/2)^2\ge 1-\varepsilon$. 

By Fubini's theorem, there exists some $c$ such that $f_c$ is domatic at a $\mu$-measure $\ge 1-\varepsilon$ set of $x\in X$. This Borel function $f_c:X\to \{0,\ldots,k-1\}$ is as desired. {}
\end{proof}

\begin{theorem}\label{thm:3.5}
Let $(X,\mu)$ be a Borel probability space, and let $G$ be an out-degree $\aleph_0$-regular Borel graph with countable in-degrees on the vertex set $X$. Then $G$ admits a $\mu$-measurable domatic $\aleph_0$-partition. 
\end{theorem}

\begin{proof}
First we explain why we may assume that $G$ is quasi-$\mu$-preserving, meaning that every $\mu$-null set is contained in an $E_G$-invariant $\mu$-null set. By the Feldman--Moore theorem \citep[\S{}4.1]{kechris-marks}, the countable Borel equivalence relation $E_G=\bigcup_{n<\omega} T_n\subseteq X^2$ is generated by countably many Borel involutions $T_n:X\to X$. Letting $\nu=\sum_{n<\omega}2^{-n-1}(T_n)_*(\mu)$, we get a Borel probability measure $\nu$ such that every $\nu$-null set is $\mu$-null, and $G$ is quasi-$\nu$-preserving. Thus $\nu$-measurable functions are $\mu$-measurable, and we may replace $\mu$ with $\nu$ to assume without loss of generality that $G$ is quasi-$\mu$-preserving. 

By Lemma \ref{lem:3.2} and since we assumed $G$ is quasi-$\mu$-preserving, it suffices to show the existence of a $\mu$-measurable function $f:X\to C$ such that $\tabs{C}=\aleph_0$, and for a $\mu$-conull set of $x\in X$, the set $f[N_G(x)]\subseteq C$ is infinite. 

For every $n<\omega$, Lemma \ref{lem:3.4} implies there is a Borel function $f_n:X\to \{0,\ldots,2^n-1\}$ that is domatic at a set of vertices $A_n\subseteq X$ of $\mu$-measure $\mu(A_n)\ge 1-2^{-n}$. The $2^n$ color classes of $f_n$ partition $X$, and so the color class $D_n\subseteq X$ of $f_n$ with the least $\mu$-measure satisfies $\mu(D_n)\le 2^{-n}$. Since $f_n$ is domatic at $A_n$, we have $D_n\cap N_G(x)\ne \varnothing$ for all $x\in A_n$. 

Since $\mu(X\ssm A_n)\le 2^{-n}$ and $\mu(D_n)\le 2^{-n}$ for all $n<\omega$, the Borel--Cantelli lemma implies that there is a $\mu$-conull set of $x\in X$ which lies in only finitely many of the sets $X\ssm A_n$ and $D_n$. Since we assumed $G$ is quasi-$\mu$-preserving, we may find an $E_G$-invariant $\mu$-conull subset $Y\subseteq X$ of such elements $x\in X$ that lie in finitely many of the sets $X\ssm A_n$ and $D_n$. 

Since every $y\in Y$ lies in finitely many $D_n$, we can define the $\mu$-measurable function $f:Y\to[\omega]^{<\omega}$ such that for every $y\in Y$, $f(y)\subseteq\omega$ is the finite set of all $n<\omega$ for which $y\in D_n$. We claim that for every $x\in Y$, the set $f[N_G(x)]\subseteq [\omega]^{<\omega}$ is infinite, which completes our proof as $C=[\omega]^{<\omega}$ has size $\aleph_0$ as we wanted. 

Let $x\in Y$. Since $x$ lies in finitely many $X\ssm A_n$, it lies in infinitely many $A_n$, and so there are infinitely many $D_n$ such that $D_n\cap N_G(x)\ne\varnothing$ by the previous arguments. By definition of $f$, infinitely many $n<\omega$ lie in the set $\bigcup_{y\in N_G(x)}f(y)=\bigcup f[N_G(x)]\subseteq\omega$. Since $f[N_G(x)]\subseteq[\omega]^{<\omega}$ is a collection of finite sets whose union $\bigcup f[N_G(x)]\subseteq\omega$ is infinite, the collection $f[N_G(x)]$ itself is infinite, which completes the proof. {} 
\end{proof}

\begin{corollary}\label{cor:3.6}
Let $\Gamma$ be a countable group with a free Borel action on a Borel probability space $(X,\mu)$. Let $S\subseteq\Gamma$ be a countably infinite subset. Then the Schreier graph $\Sch(\Gamma,S,X)$ admits a $\mu$-measurable domatic $\aleph_0$-partition. 
\end{corollary}

\begin{proof}
Follows from Theorem \ref{thm:3.5}. {}
\end{proof}

See also Theorem \ref{thm:4.3} for the Baire category counterpart of Theorem \ref{thm:3.5}.

\subsection{Domatic edge-\texorpdfstring{$\aleph_0$}{infinite}-partitions}\label{sec:3.3}

In this section, we explore an edge-coloring version of domatic partitions. 

\begin{definition}\label{def:3.7}
Let $G$ be a loop-free simple undirected graph on a vertex set $V$, and we write $G\subseteq V^2$ for its set of edges as before. Let $f:G\partialto{}C$ be a symmetric partial function defined on $G$'s set of edges, meaning that $f(w,v)=f(v,w)$ for every $(v,w)\in\dom(f)$. We say that $f$ is \emph{domatic} at a vertex $v\in V$ if for every $c\in C$ there exists $w\in N_G(v)$ such that $f(v,w)=f(w,v)=c$. The symmetric partial function $f:G\partialto C$ is \emph{domatic} if it's domatic everywhere in $V$. 
\end{definition}

Let $G$ be a loop-free simple undirected graph on a vertex set $V$, and let $G\subseteq V^2$ be its set of edges. Let $\sim$ be the equivalence relation on $G$ generated by $(v,w)\sim(w,v)$, so that $G/{\sim}$ is the set of undirected edges of $G$. We may define the \emph{subdivision} of $G$ as the new graph $G'$, such that $G'$ is bipartite on the vertex set $V\sqcup (G/{\sim})$, and the edges of $G'$ are pairs $(v,[e]_{\sim})$ where $v\in V$, $e\in G$, and $e$ is incident to $v$. Note that if $G$ is Borel then $G'$ is Borel and Borel bipartite. 

We see that a symmetric partial function $f:G\partialto C$ is domatic as an edge-coloring for $G$, if and only if the quotient $f/{\sim}:(G/{\sim})\partialto{}C$ is domatic at $V$ as a vertex-coloring for $G'$. Also when $G$ is locally finite or locally countable, $G'$ is also locally finite or locally countable respectively. Thus we can transfer results about domatic vertex-partitions to results about domatic edge-partitions, using the auxiliary graph $G'$. 

\begin{theorem}\label{thm:3.8}
Let $G$ be a loop-free simple undirected $\aleph_0$-regular Borel graph on a Borel space $X$ of vertices. If $\mu$ is any Borel probability measure on $X$, then there is an $E_G$-invariant $\mu$-conull Borel set $C_\mu\subseteq X$ and a symmetric Borel function $f_\mu:G\restrict{}C_\mu\to \omega$ such that $f_\mu$ is domatic everywhere in $C_\mu$. Similarly, if $\tau$ is any Polish topology on $X$, then there is an $E_G$-invariant $\tau$-comeager Borel set $C_\tau\subseteq X$ and a symmetric Borel function $f_\tau:G\restrict{} C_\tau\to \omega$ such that $f_\tau$ is domatic everywhere in $C_\tau$. 
\end{theorem}

\begin{proof}
By Feldman--Moore \citep[Proposition 4.1]{kechris-marks}, $G$ admits a symmetric Borel edge-$\aleph_0$-coloring $f:G\to\omega$, meaning that for every two distinct edges $e,e'\in G$ that share a common vertex, we have $f(e)\ne f(e')$. In particular, every vertex of $G$ belongs to edges of infinitely many $f$-colors. Then the proofs of Lemmas \ref{lem:3.2} and \ref{lem:3.3} imply our desired result. {}
\end{proof}

We will continue our discussions about domatic edge-partitions in Section \ref{sec:4.3}, where we prove Theorem \ref{thm:4.5} as a Borel counterpart of Theorem \ref{thm:3.8}. 

\section{More examples of domatic finite partitions}\label{sec:4}

\subsection{Maximal independent sets are domatic 2-partitions}\label{sec:4.1}

In Section \ref{sec:1.1}, we defined a subset $I\subseteq V$ to be \emph{independent} for a graph with loops $G\subseteq V^2$ if for every edge $(v,w)\in G$ with $v\ne w$, not both $v$ and $w$ belong to $I$. Given a graph with loops $G\subseteq V^2$, we say that a vertex $v\in V$ is \emph{isolated} if $N_G(v)\subseteq\{v\}$. 

\begin{lemma}\label{lem:4.1}
Let $G$ be a fully looped undirected graph without isolated vertices on a vertex set $V$. If $A\subseteq V$ is a maximal $G$-independent set, then $A$ and $V\ssm A$ are dominating sets for $G$ which form a domatic $2$-partition. 
\end{lemma}

\begin{proof}
First we show that $A$ is dominating. Let $v\in V$ and it suffices to show $A\cap N_G(v)\ne\varnothing$. If not, then $A\cap N_G(v)=\varnothing$, and so $A\cup\{v\}\supsetneqq A$ is also a $G$-independent set, which contradicts maximality of $A$.

Next we show that $V\ssm A$ is dominating. Let $v\in V$ and it suffices to show $(V\ssm A)\cap N_G(v)\ne\varnothing$. Since $v$ is not isolated in $G$, it has a neighbor $w\in N_G(v)\ssm\{v\}$. Since $A$ is independent, either $v$ or $w$ belongs to $(V\ssm A)\cap N_G(v)\ne \varnothing$ as desired. {}
\end{proof}

\begin{theorem}\label{thm:4.2}
Let $G$ be a fully looped undirected locally countable Borel graph without isolated vertices on a Borel space $X$ of vertices. Then $G$ admits a Borel domatic $2$-partition. 
\end{theorem}

\begin{proof}
By Lusin--Novikov uniformization \citep[Theorem 18.10]{kechris}, there is a Borel function $f:X\to X$ such that for every $x\in X$, $f(x)\in N_G(x)\ssm\{x\}$ is a neighbor of $x$. Let $G_f$ be the fully looped undirected Borel subgraph of $G$ generated by loops and edges of the form $(x,f(x))$ for $x\in X$. By Kechris--Solecki--Todorcevic \citep[Corollary 4.6]{kechris-marks}, $G_f$ has a Borel $\aleph_0$-coloring, and by Kechris--Solecki--Todorcevic \citep[Proposition 4.9]{kechris-marks}, $G_f$ has a Borel maximal independent set. By Lemma \ref{lem:4.1}, $G_f$ has a Borel domatic $2$-partition, and since $G$ has $G_f$ as a subgraph, $G$ also has the same Borel domatic $2$-partition. 
\end{proof}

\subsection{Baire measurable domatic 3-partitions}\label{sec:4.2}

The following result is a Baire category counterpart of Theorem \ref{thm:3.5}. 

\begin{theorem}\label{thm:4.3}
There exists a fully looped undirected $\aleph_0$-regular acyclic Borel graph $G$ on a Polish space $(X,\tau)$ of vertices, without $\tau$-Baire measurable domatic $3$-partitions. 
\end{theorem}

\begin{proof}
The graph $G$ is Lecomte's infinite dimensional version of the Kechris--Solecki--Todorcevic graph $G_0$ \citep{miller11}. We give a full proof below for sake of self-containedness. 

Let $\langle s_n:n<\omega\rangle$ be a sequence of finite sequences $s_n\in {}^{\omega>}\omega$, such that $\tabs{s_n}=n$ for all $n<\omega$, and $\{s_n:n<\omega\}$ is dense in ${}^{\omega>}\omega$. Let $\{N_s:s\in {}^{\omega>}\omega\}$ be the standard topological basis of ${}^\omega \omega$. The graph $G$ will be bipartite on the vertex set $X={}^{\omega}\omega\sqcup\bigsqcup_{n<\omega}(N_{s_n\concat{}0}\times\{s_n\})$, where we equip $X$ with its natural Polish topology $\tau$. Non-loop edges of $G$ are generated by pairs $(s_n\concat{}k\concat{}x,(s_n\concat{}0\concat{}x,s_n))$ such that $k,n<\omega$ and $x\in {}^\omega \omega$. By usual arguments, one sees that the $E_G$-saturation of $\tau$-meager sets are $\tau$-meager sets. 

We'll check that there is an $E_G$-invariant $\tau$-comeager $G_\delta$ set $C\subseteq X$, such that the graph $G\restrict{} C$ is $\aleph_0$-regular. Each vertex in the second part $\bigsqcup_n (N_{s_n\concat{}0}\times\{s_n\})$ always has $G$-degree $\aleph_0$, whereas a vertex $x\in {}^\omega \omega$ in the first part of $X$ has $G$-degree $\aleph_0$ if and only if $s_n\preceq x$ for infinitely many $n<\omega$. Since $\{s_n:n<\omega\}$ is dense in ${}^{\omega>}\omega$, there is a $\tau$-comeager $G_\delta$ set of $x\in {}^\omega \omega\subseteq X$ which has $G$-degree $\aleph_0$, from which we get our desired $E_G$-invariant $\tau$-comeager $G_\delta$ set $C\subseteq X$ using that meager sets saturate to meager sets. 

The proof that $G$ is acyclic is the same as the argument that the Kechris--Solecki--Todorcevic graph $G_0$ is acyclic \citep[Example 4.16]{kechris-marks}, which we omit here. 

Finally, we'll check that for any $E_G$-invariant $\tau$-comeager set $C\subseteq X$, the graph $G\restrict{}C$ admits no $\tau$-Baire measurable domatic $3$-partitions. Assume that $f:C\to \{0,1,2\}$ is a $\tau$-Baire measurable domatic function. Then by density of $\{s_n:n<\omega\}\subseteq {}^{\omega>}\omega$, there exists some basic open set $N_{s_n}\subseteq {}^\omega\omega$ and an $E_G\restrict{} N_{s_n}$-invariant $\tau$-comeager subset $B\subseteq N_{s_n}\cap C$ of $N_{s_n}$, such that $f$ is constant over $B$. Let $z=s_n\concat{}0\concat{}x\in B$, and we note that $N_G(z,s_n)\subseteq B\cup\{(z,s_n)\}$ by $E_G\restrict{}N_{s_n}$-invariance of $B$. So $f$ can take on at most two output values over $N_G(z,s_n)$, and in particular $f$ is not domatic at $(z,s_n)\in C$. {}
\end{proof}

\subsection{Borel domatic edge-2-partitions}\label{sec:4.3}

We now resume our discussions about domatic edge-partitions from Section \ref{sec:3.3}. The result Theorem \ref{thm:4.5} proved by Felix Weilacher is a Borel counterpart of Theorem \ref{thm:3.8}. 

\begin{lemma}[Weilacher]\label{lem:4.4}
Let $\F_2=\langle a,b\rangle$ be the free group on $2$ generators $a,b$. There exists a free Borel action of $\F_2$ on a Borel space $X$, such that the countable Borel equivalence relations $E_{\langle a\rangle}^X,E_{\langle b\rangle}^X$ are smooth, and every Borel function $f:X\to\{0,1\}$ admits either a $0$-monochromatic $\langle a\rangle$-orbit or a $1$-monochromatic $\langle b\rangle$-orbit. 
\end{lemma}

\begin{proof}
The proof strategy is to modify Marks' Borel determinacy lemma \citep[Lemma 2.1]{marks}. Note that by Marks \citep[Lemma 2.1]{marks}, the free Borel $\F_2$-space $\operatorname{Free}(\omega^{\F_2})$ already satisfies that every Borel function $f:\operatorname{Free}(\omega^{\F_2})\to\{0,1\}$ has a $0$-monochromatic $\langle a\rangle$-orbit or a $1$-monochromatic $\langle b\rangle$-orbit, and it suffices to modify the proof to guarantee also the smoothness of $E_{\langle a\rangle}^X,E_{\langle b\rangle}^X$. 

Let $X\subseteq\operatorname{Free}(\omega^{\F_2})$ be the Borel subset of all $x\in \operatorname{Free}(\omega^{\F_2})$ such that for every $\gamma\in \F_2$, the functions $(\gamma\cdot x)\restrict{}\langle a\rangle:\langle a\rangle\to\omega$ and $(\gamma\cdot x)\restrict{}\langle b\rangle:\langle b\rangle\to\omega$ are injective, that is, $X$ is the set of functions $\F_2\to\omega$ in $\operatorname{Free}(\omega^{\F_2})$ which are injective over every coset of $\langle a\rangle$ or $\langle b\rangle$. The proof of Marks' lemma \citep[Lemma 2.1]{marks} applies to $X$, since there we may require that both players of the game make moves that are partial functions of functions in $X$. It follows that every Borel function $f:X\to\{0,1\}$ has a $0$-monochromatic $\langle a\rangle$-orbit or a $1$-monochromatic $\langle b\rangle$-orbit. 

Next, note that the countable Borel equivalence relations $E_{\langle a\rangle}^X,E_{\langle b\rangle}^X$ admit Borel selectors, since in each equivalence class we can select the unique $x:\F_2\to\omega$ which minimizes the value $x(1_{\F_2})<\omega$ by definition of $X$. Thus $E_{\langle a\rangle}^X,E_{\langle b\rangle}^X$ are smooth. {}
\end{proof}

In Section \ref{sec:3.3}, we defined a \emph{symmetric} function $f:G\to C$ on a graph $G\subseteq V^2$ to be such that $f(v,w)=f(w,v)$ for every $(v,w)\in G$, and we also defined the equivalence relation $\sim$ on a loop-free simple undirected graph $G\subseteq V^2$ to be generated by $(v,w)\sim (w,v)$.  

\begin{theorem}[Weilacher]\label{thm:4.5}
There exists a loop-free simple undirected $\aleph_0$-regular acyclic Borel graph $G\subseteq X^2$ that is Borel bipartite on a vertex set $X=A\sqcup B$, such that every symmetric Borel function $f:G\to\{0,1\}$ admits either a vertex $a\in A$ belonging to only edges of color $0$, or a vertex $b\in B$ belonging to only edges of color $1$. In particular, $f$ is not domatic at such a vertex $a$ or $b$. 
\end{theorem}

\begin{proof}
Let $X_0$ be the Borel $\F_2$-space given by Lemma \ref{lem:4.4}. Since $E_{\langle a\rangle}^{X_0},E_{\langle b\rangle}^{X_0}$ are smooth, the quotient spaces $A=X_0/E_{\langle a\rangle}^{X_0}$ and $B=X_0/E_{\langle b\rangle}^{X_0}$ are Borel. We will define the graph $G$ over the vertex set $A\sqcup B$, such that the edges of $G$ are generated by pairs $([x]_{E_{\langle a\rangle}^{X_0}},[x]_{E_{\langle b\rangle}^{X_0}})$ for some $x\in X_0$. Since the $\F_2$-action on $X_0$ is free, we see that $G$ is $\aleph_0$-regular acyclic. 

Note that the function $x\mapsto [([x]_{E_{\langle a\rangle}^{X_0}},[x]_{E_{\langle b\rangle}^{X_0}})]_{\sim}$ defines a Borel isomorphism $X_0\cong G/{\sim}$, under which we may view symmetric Borel functions $f:G\to\{0,1\}$ as Borel functions $F:X_0\to\{0,1\}$. Lemma \ref{lem:4.4} implies that such a Borel function $F$ is either $0$-monochromatic on an $E_{\langle a\rangle}^{X_0}$-class $a_0\in A$, or it is $1$-monochromatic on an $E_{\langle b\rangle}^{X_0}$-class $b_0\in B$. This means that either $a_0\in A$ belongs to only $G$-edges of $f$-color $0$, or $b_0\in B$ belongs to only $G$-edges of $f$-color $1$, as we desired. {}
\end{proof}

It's not hard to see, via a chase of equivalence between definitions, that the graph $G$ from Theorem \ref{thm:4.5} is exactly one which admits no Borel sinkless orientations. 

\subsection{Locally finite graphs}\label{sec:4.4}

In this section, we're finally able to move our attention away from $\aleph_0$-regular Borel graphs. The main challenge we face when constructing measurable domatic partitions for locally finite Borel graphs is a lack of good results on when even finite graphs admit domatic partitions. Nevertheless, the following are a few selected examples among things one could say about domatic partitions for fully looped undirected locally finite Borel graphs:

\begin{enumerate}[(1)]
\item 
Regular graphs of sufficiently large finite degree. 

It's a standard application of the Lov\'{a}sz local lemma in probabilistic combinatorics that for a fixed $k<\omega$ and every sufficiently large $d<\omega$, every fully looped undirected $d$-regular graph admits domatic $k$-partitions. Results by Bernshteyn \citep[Theorem 2.20]{bernshteyn} and Cs\'{o}ka--Grabowski--M\'{a}th\'{e}--Pikhurko--Tyros \citep[Theorem 4.5]{cgmpt} imply that versions of the Lov\'{a}sz local lemma still hold for various classes of coloring problems in the context of descriptive graph combinatorics. It follows more or less directly that for a fixed $k<\omega$ and every sufficiently large $d<\omega$, every fully looped undirected $d$-regular Borel graph admits measure-theoretic and Baire measurable domatic $k$-partitions unconditionally, and Borel domatic $k$-partitions as long as the graph is of uniform subexponential growth. 

The above technique of using the Lov\'{a}sz local lemma can also be applied to find domatic finite partitions for $\aleph_0$-regular Schreier graphs. Let $\Gamma$ be a countably infinite group, $S\subseteq \Gamma$ a countably infinite generating subset, and $X$ a free Borel $\Gamma$-space, so that the directed Schreier graph $G=\Sch(\Gamma,S,X)$ is an out-degree $\aleph_0$-regular, in-degree $\aleph_0$-regular Borel graph on $X$. For every $k<\omega$, if $F\subseteq S$ is a sufficiently large finite set, then $\Sch(\Gamma,F,X)$ is a locally finite regular subgraph of $G$ of sufficiently large finite degree. So the measurable Lov\'{a}sz local lemmas imply that $\Sch(\Gamma,F,X)$ and hence $G=\Sch(\Gamma,S,X)$ admit measure-theoretic and Baire measurable domatic $k$-partitions. Moreover, $G=\Sch(\Gamma,S,X)$ admits Borel domatic $k$-partitions for all $k<\omega$ if every finitely generated subgroup of $\Gamma$ has subexponential growth, and we don't know yet if there exists an example of a free Borel $\Gamma$-space $X$ for which the $\aleph_0$-regular graph $\Sch(\Gamma,S,X)$ does not admit a Borel domatic $k$-partition for some $k<\omega$. 

\item 
Locally finite acyclic graphs. 

Recall from Section \ref{sec:1.1} that if $G\subseteq V^2$ is a fully looped undirected graph on a vertex set $V$ and $G\ssm \Delta_V$ is its loop-free version, then a vertex $v\in V$ has $G$-degree $d+1$ if and only if it has $(G\ssm\Delta_V)$-degree $d$. 

The problem of finding domatic partitions on fully looped undirected locally finite acyclic graphs can be solved by the method of path decompositions by Conley--Marks--Unger \citep[Definition 1.4]{conley-marks-unger}. Given a fully looped undirected locally finite acyclic graph with a path decomposition into sufficiently long paths, one can build a straightforward greedy algorithm on each single path prioritizing its two endpoints, so that the greedy algorithm always outputs domatic coloring functions which waste at most one extra color at every vertex. It follows from the existence of sufficiently long path decompositions \citep[Lemma 3.4]{conley-marks-unger} that if $\delta(G)$ is the minimum degree of a fully looped undirected locally finite acyclic Borel graph $G$, then $G$ admits Baire measurable domatic $\delta(G)$-partitions whenever $\delta(G)\ne 3$. 

When $\delta(G)=3$, rigidity of domatic $3$-partitions for fully looped undirected bi-infinite paths implies that these domatic $3$-partitions are $3$-periodic colorings, and hence it's easy to construct free Polish $\Z$-spaces $X$ on which $\Sch(\Z,\{-1,0,1\},X)$ admits no Baire measurable domatic $3$-partitions. For example, $X$ can be the dyadic odometer. On the other hand, this case $\delta(G)=3$ always admits Borel domatic $2$-partitions by Theorem \ref{thm:4.2}. 

When we additionally assume the maximum degree $\Delta(G)$ of $G$ is bounded, the same problem of finding Baire measurable domatic $\delta(G)$-partitions for $G$ when $\delta(G)\ne 3$ can also be solved by a \textsf{TOAST} algorithm. See for example the article \citep{bcggrv} for a reference on \textsf{TOAST} algorithms. 

In the case of Borel domatic partitions for fully looped undirected locally finite acyclic graphs, one can use the Borel determinacy approach by Marks in a similar way as in Theorem \ref{thm:4.5}. For example, it follows from the analysis of $\operatorname{Free}(\omega^{(\Z/n\Z)^{*n}})$ that there is a fully looped undirected $(n+1)$-regular acyclic Borel graph without Borel domatic $3$-partitions. 
\end{enumerate}

\end{document}